\documentclass[10pt]{article}
\pdfoutput=1

\usepackage{amsmath}
\usepackage{mathtools}      
\usepackage{stackrel}
\usepackage{scalerel}
\usepackage{leftidx}
\usepackage{xfrac}

\usepackage[charter]{mathdesign}   
\usepackage{euscript}              

\let\mathbb\undefined
\usepackage{bbold}                 

\DeclareSymbolFont{usualmathcal}{OMS}{cmsy}{m}{n}
\DeclareSymbolFontAlphabet{\mathcal}{usualmathcal}

\usepackage{geometry}
\usepackage{xcolor}                
\usepackage{tcolorbox}
\usepackage{graphicx}
\usepackage{subfigure}             
\usepackage{verbatim}
\usepackage{comment}

\usepackage{authblk}

\usepackage[shortlabels]{enumitem}
\usepackage{appendix}

\usepackage[all,2cell,arrow,matrix]{xy}
\UseAllTwocells
\SilentMatrices
\usepackage{tikz-cd}
\usepackage{tikz}
\usetikzlibrary{arrows,arrows.meta,decorations.markings,decorations.pathreplacing,decorations.pathmorphing,calc}

\tikzset{->-/.style={decoration={markings,mark=at position #1 with {\arrow{Stealth}}},postaction={decorate}},->-/.default=0.55}

\usepackage[nottoc]{tocbibind}     
\usepackage{hyperref}
\usepackage{cleveref}              
\hypersetup{colorlinks=true, linkcolor=red!50!black, citecolor=green!50!black, urlcolor=blue!80!black}

\usepackage[amsmath,amsthm,thmmarks]{ntheorem}
\theoremstyle{definition}

\newtheorem{mainthm}{Theorem}

\newtheorem{thm}{Theorem}[section]
\newtheorem{prop}[thm]{Proposition}

\newtheorem{cor}[thm]{Corollary}
\newtheorem{lem}[thm]{Lemma}

\newtheorem{conj}[thm]{Conjecture}

{
\theoremsymbol{\mbox{${\blacksquare}$}}
\newtheorem{defn}[thm]{Definition}
}
{
\theoremsymbol{\mbox{$\heartsuit$}}
\newtheorem{expl}[thm]{Example}
}
{
\theoremsymbol{\mbox{$\diamondsuit$}}
\newtheorem{rem}[thm]{Remark}
}
\qedsymbol{\mbox{$\square$}}

\numberwithin{equation}{section}
\numberwithin{thm}{section}

\newcommand{\be}{\begin{equation}}
\newcommand{\ee}{\end{equation}}
\newcommand{\bea}{\begin{eqnarray}}
\newcommand{\eea}{\end{eqnarray}}
\newcommand{\bnu}{\begin{enumerate}}
\newcommand{\enu}{\end{enumerate}}
\newcommand{\bit}{\begin{itemize}}
\newcommand{\eit}{\end{itemize}}

\newcommand{\pf}{\begin{proof}}
\newcommand{\epf}{\end{proof}}

\makeatletter
\providecommand{\leftsquigarrow}{%
  \mathrel{\mathpalette\reflect@squig\relax}%
}
\newcommand{\reflect@squig}[2]{%
  \reflectbox{$\m@th#1\rightsquigarrow$}%
}
\makeatother

\newcommand{\Cb}{\mathbb{C}}
\newcommand{\Hb}{\mathbb{H}}

\newcommand{\Rb}{\mathbb{R}}
\newcommand{\Zb}{\mathbb{Z}}
\newcommand{\bk}{\mathbb{k}}
\newcommand{\Eb}{\mathbb{E}}
\newcommand{\Fb}{\mathbb{F}}

\newcommand{\CB}{\EuScript{B}}
\newcommand{\CC}{\EuScript{C}}
\newcommand{\CD}{\EuScript{D}}
\newcommand{\CE}{\EuScript{E}}

\newcommand{\CG}{\EuScript{G}}

\newcommand{\CM}{\EuScript{M}}
\newcommand{\CN}{\EuScript{N}}

\newcommand{\CP}{\EuScript{P}}

\newcommand{\CX}{\EuScript{X}}

\newcommand{\SF}{\mathsf{F}}

\newcommand{\rmB}{\mathrm{B}}

\newcommand{\rmM}{\mathrm{M}}

\DeclareMathOperator{\Hom}{Hom}
\DeclareMathOperator{\End}{End}
\DeclareMathOperator{\Aut}{Aut}

\DeclareMathOperator{\ev}{ev}

\DeclareMathOperator{\ob}{ob}

\DeclareMathOperator{\fun}{Fun}

\DeclareMathOperator{\Alg}{Alg}

\DeclareMathOperator{\LMod}{LMod}
\DeclareMathOperator{\RMod}{RMod}
\DeclareMathOperator{\BMod}{BMod}

\DeclareMathOperator{\Gal}{Gal}

\newcommand{\op}{\mathrm{op}}
\newcommand{\rev}{\mathrm{rev}}

\newcommand{\one}{\mathbb{1}}

\newcommand{\forget}{\text{\usefont{U}{euf}{m}{n}f}}

\newcommand{\vect}{\mathrm{Vec}}
\newcommand{\svect}{\mathrm{sVec}}

\newcommand{\rep}{\mathrm{Rep}}
\newcommand{\srep}{\mathrm{sRep}}

\newcommand{\bscale}{0.7}
\makeatletter
\newcommand{\ec}[2][]{{\@ec{#1 |}{#2}}}
\newcommand{\bc}[2][]{{\@ec{#1}{#2}}}
\newcommand{\@ec}[2]{\mathchoice
  {\displaystyle \raise.9ex\hbox{$\scaleobj{\bscale}{#1}$} {#2}}%
  {\textstyle \raise.9ex\hbox{$\scaleobj{\bscale}{#1}$} {#2}}%
  {\scriptstyle \raise.55ex\hbox{$\scriptstyle \scaleobj{\bscale}{#1}$} {#2}}%
  {\scriptscriptstyle \raise.38ex\hbox{$\scriptscriptstyle \scaleobj{\bscale}{#1}$} {#2}}%
}
\makeatother

\newcommand{\dquotient}{/\mkern-6mu/}
\newcommand{\rex}{\mathrm{rex}}
\newcommand{\sep}{\mathrm{sep}}
\newcommand{\rig}{\mathrm{rig}}

\title{\huge Classification of symmetric fusion categories over \texorpdfstring{$\Rb$}{R}}

\author[a,b]{Mo Huang \thanks{Email: \href{mailto:jasmine.huang.527@gmail.com}{\tt jasmine.huang.527@gmail.com}}}
\author[c]{Hao Xu \thanks{Email: \href{mailto:haoxu@imada.sdu.dk}{\tt haoxu@imada.sdu.dk}}}
\author[b]{Zhi-Hao Zhang \thanks{Email: \href{mailto:zhangzhihao@bimsa.cn}{\tt zhangzhihao@bimsa.cn}}}
\affil[a]{School of Mathematical Sciences, East China Normal University, Shanghai, 200241, China}
\affil[b]{Beijing Institute of Mathematical Sciences and Applications, Beijing, 101408, China}
\affil[c]{Centre for Quantum Mathematics, University of Southern Denmark, Odense, 5230, Denmark}

\date{\vspace{-5ex}}

\begin{document}

\maketitle

\begin{abstract}
We show that every symmetric fusion category over $\Rb$ is equivalent to the category of finite-dimensional semi-linear representations of a $\Zb_2$-graded finite super group. The proof uses Galois descent for tensor categories over $\Cb/\Rb$, reducing the classification to semi-linear $\Zb_2$-actions on symmetric fusion categories over $\Cb$. As a further structural result, we establish a Tannaka-Krein type correspondence between symmetric fusion categories over $\Rb$ and finite groupoids with a $\Zb_2 \times \rmB \Zb_2$-action. This gives a complete real analogue of Deligne's classification result.
\end{abstract}

\tableofcontents

\section{Introduction}

In quantum mechanics, symmetry plays a fundamental role. Typically, the state space of a quantum system is a Hilbert space $\mathcal{H}$, and a symmetry is realized as a unitary action of a group $G$ on $\mathcal{H}$. When $G$ is a finite group (or a compact Lie group), all its irreducible complex representations are finite-dimensional, and $\mathcal{H}$ decomposes into a direct sum of irreducible subrepresentations. The representation category $\rep_\Cb(G)$ thus becomes a central object for studying such symmetries. 

For fermionic systems, the state space $\mathcal{H}$ carries a natural $\Zb_2$-grading $\mathcal{H} = \mathcal{H}_0 \oplus \mathcal{H}_1$, where $\mathcal{H}_0$ and $\mathcal{H}_1$ consist of states with even and odd number of fermions, respectively. There is a fermion parity operator $(-1)^f$ acting as $1$ on $\mathcal{H}_0$ and $-1$ on $\mathcal{H}_1$. A symmetry group $G$ acting on $\mathcal{H}$ should preserve the parity of the fermion number, i.e., it must commute with $(-1)^f$. Consequently, $G$ and the fermion parity generate a central extension $(\tilde{G},z)$, where $z \colon \Zb_2 \to \tilde{G}$ is an embedding, or equivalently, a central element $z \in \tilde{G}$ of order two. The original symmetry then corresponds to a representation of $\tilde{G}$ on $\mathcal{H}$ such that $z$ acts as the fermion parity. We call $(\tilde{G}, z)$ a \emph{super group} and such a representation a \emph{super representation}. The category of finite-dimensional super representations is denoted by $\srep_\Cb(\tilde{G},z)$. Under the tensor product $\otimes_\Cb$ of $\Cb$-vector spaces, $\srep_\Cb(\tilde{G},z)$ is a tensor category equivalent to $\rep_\Cb(\tilde{G})$, but it has a different braiding structure because exchanging two fermions incurs a phase of $-1$. 

When the groups are finite, these categories are examples of symmetric fusion categories. Deligne's theorem \cite{Del90,Del02} tells us that every symmetric fusion category over $\Cb$ is equivalent to either $\rep_\Cb(G)$ for some finite group $G$, or $\srep_\Cb(G,z)$ for some finite super group $(G,z)$. Symmetric fusion categories therefore abstractly capture the notion of finite unitary symmetries.

Wigner's theorem, however, asserts that symmetries of a quantum system can be unitary or anti-unitary. For example, time-reversal symmetry is anti-unitary. Hence a symmetry group $G$ generally decomposes as $G = G_0 \sqcup G_1$, where elements of $G_0$ act unitarily and those of $G_1$ act anti-unitarily. Equivalently, this decomposition is given by a group homomorphism $s \colon G \to \Zb_2$ with $G_0 = \ker(s)$. We call $(G,s)$ a \emph{$\Zb_2$-graded group} and the corresponding representations \emph{semi-linear} or \emph{semi-unitary representations}. The category of finite-dimensional semi-linear representations of $(G,s)$ is denoted by $\rep_{\Cb/\Rb}(G,s)$. When $G$ is finite, this is again a symmetric fusion category, but it is only $\Rb$-linear in general. Similarly, for a general symmetry on a fermionic system, the symmetry group is a \emph{$\Zb_2$-graded super group} $(G,z,s)$, which is both a super group and a $\Zb_2$-graded group such that $z \in \ker(s)$ (because the fermion parity operator is always unitary). The category of finite-dimensional \emph{semi-linear super representations} is denoted by $\srep_{\Cb/\Rb}(G,z,s)$. The main result of this paper is that all symmetric fusion categories over the real numbers $\Rb$ are of these forms. Therefore, an abstract symmetry of a quantum system can be described by a symmetric fusion category over $\Rb$.

\begin{mainthm}[Theorem \ref{thm_classification_symmetric_fusion_R}]
Every symmetric fusion category over $\Rb$ is equivalent to $\srep_{\Cb/\Rb}(G,z,s)$ for some $\Zb_2$-graded finite super group $(G,z,s)$.
\end{mainthm}

The proof proceeds via Galois descent for categories over the extension $\Cb/\Rb$. Given a Galois extension $\mathbb{E}/\mathbb{F}$, Galois descent is the standard way to study algebraic structures over the smaller field $\mathbb{F}$ from those over the larger field $\mathbb{E}$: extending scalars turns an $\mathbb{F}$-structure into an $\mathbb{E}$-structure, while an $\mathbb{E}$-structure equipped with suitable descent data corresponds to an $\mathbb{F}$-structure. This is classical textbook material for structures such as vector spaces or algebras. The descent theory for (tensor) categories has been developed in \cite{EG11}. In the present paper we carefully analyze the descent data for categories over the extension $\Cb/\Rb$ and show that it is given by a semi-linear action of $\Zb_2 \simeq \Gal(\Cb/\Rb)$. The results can be summarized as a symmetric monoidal equivalence of 2-categories.

\begin{mainthm}[Theorem \ref{thm_Galois_descent_category}]
Complexification and equivariantization define an $\Rb$-linear symmetric monoidal equivalence $2\vect_\Cb^{\Zb_2} \simeq 2\vect_\Rb$, where the $\Zb_2$-action on $2\vect_\Cb$ is defined by complex conjugation.
\end{mainthm}

As a consequence, a symmetric fusion category $\CC$ over $\Rb$ is equivalent to the equivariantization of its complexification, which is a symmetric fusion category over $\Cb$ if $\CC(\one,\one) \simeq \Rb$. Therefore, by Deligne's theorem \cite{Del90,Del02}, it suffices to classify semi-linear $\Zb_2$-actions on $\srep_\Cb(K,z)$ for finite super groups $(K,z)$. We explicitly compute these actions and show that they correspond to $\Zb_2$-graded extensions of $(K,z)$, thereby proving our main classification result.

Besides the concrete classification, we also establish a Tannaka-Krein type reconstruction theorem for symmetric (multi-)fusion categories over $\Rb$. Given a symmetric fusion category $\CC$ over $\Rb$, a \emph{fiber functor} is a faithful exact symmetric monoidal $\Rb$-linear functor from $\CC$ to $\vect_\Cb$ or $\svect_\Cb$. In physics, fiber functors correspond to concrete Hilbert space realizations of the abstract symmetry category. All $\Rb$-linear fiber functors $\CC \to \svect_\Cb$ form a groupoid. It admits an action of the 2-group of $\Rb$-linear autoequivalences of $\svect_\Cb$, which is equivalent to $\Zb_2 \times \rmB \Zb_2$. Conversely, given a finite groupoid $X$ equipped with a $\Zb_2 \times \rmB \Zb_2$-action, the category $\fun_{\Zb_2 \times \rmB \Zb_2}(X,\svect_\Cb)$ of equivariant functors is a symmetric multi-fusion category over $\Rb$. These two constructions are mutually inverse.

\begin{mainthm}[Theorem \ref{thm_symmetric_fusion_category_R_super_groupoid_Galois_action}]
Taking the groupoid of $\Rb$-linear fiber functors defines an opposite equivalence between the 2-groupoids of symmetric multi-fusion categories over $\Rb$ and the 2-groupoids of finite groupoids equipped with a $\Zb_2 \times \rmB \Zb_2$-action.
\end{mainthm}

As a byproduct of the reconstruction theorem, we also characterize exactly which symmetric fusion categories over $\Rb$ admit fiber functors to $\vect_\Rb$ (Tannakian) or to $\svect_\Rb$ (super Tannakian).

Many of the discussions in this paper work equally well for arbitrary finite Galois extensions, as in \cite{EG11}. For infinite Galois extensions, the Galois group is a profinite group and one needs to consider smooth (continuous) actions on categories, which makes the situation substantially more complicated. We will treat the general Galois descent of tensor categories and the classification of symmetric fusion categories over fields of characteristic zero in a subsequent paper.

\medskip
The organization of this paper is as follows. In Section \ref{sec_preliminary} we review some basic results of fusion categories over general fields. In Section \ref{sec_Deligne_symmetric_fusion}, we review Deligne's classification of symmetric fusion categories over an algebraically closed field of characteristic zero. We also establish an opposite equivalence between the 2-groupoids of symmetric multi-fusion categories and the 2-groupoids of finite super groupoids. In Section \ref{sec_Galois_descent}, we discuss the Galois descent for categories over the extension $\Cb/\Rb$ and summarize the results into a symmetric monoidal equivalence of 2-categories. In Section \ref{sec_semi-linear_rep}, we study the representation categories of $\Zb_2$-graded (super) groups, which are examples of symmetric fusion categories over $\Rb$. Finally, in Section \ref{sec_symmetric_fusion_R}, we classify symmetric fusion categories over $\Rb$ by computing semi-linear $\Zb_2$-actions on symmetric fusion categories over $\Cb$. We also determine $\Rb$-linear fiber functors to $\svect_\Cb$, and characterize symmetric fusion categories admitting fiber functors to $\vect_\Rb$ or $\svect_\Rb$.

\medskip
\noindent \textbf{Acknowledgements}: HX is supported by Villum Fonden 00060714 ``Global Categorical Symmetries and Phases of Quantum Matter''. ZHZ is supported by the start-up grant of BIMSA and the China Postdoctoral Science Foundation under Grant Number 2025M783072.

\section{Preliminaries} \label{sec_preliminary}

\subsection{Fusion categories}

In this subsection, we review the definition and basic results of fusion categories over a general field $\bk$. We refer readers to \cite{San25} for more details.

\begin{defn}
A $\bk$-linear abelian category is \emph{finite} if it is equivalent to the category $\LMod_A(\vect_\bk)$ of finite-dimensional left modules over some finite-dimensional $\bk$-algebra $A$. In addition, it is \emph{semisimple} or \emph{separable} if the $\bk$-algebra $A$ is semisimple or separable, respectively.
\end{defn}

In particular, a finite separable category is finite semisimple. If $\bk$ is a perfect field, every finite semisimple category is also finite separable.

\begin{rem}
The semisimplicity of a finite $\bk$-linear category is independent of the base field $\bk$, but the separability depends on the choice of $\bk$.
\end{rem}

Let $\CC,\CD,\CE$ be finite categories over $\bk$. A functor $\CC \times \CD \to \CE$ is called \emph{right exact $\bk$-linear} if it is right exact and $\bk$-linear in each variable. The category of right exact $\bk$-linear functors $\CC \times \CD \to \CE$ is denoted by $\fun_\bk^\rex(\CC,\CD;\CE)$.

\begin{defn}
The \emph{Deligne tensor product} $\CC \boxtimes_\bk \CD$ (or simply $\CC \boxtimes \CD$) is a finite category defined by the following universal property:
\bit
\item It is equipped with a right exact $\bk$-linear functor $\boxtimes \colon \CC \times \CD \to \CC \boxtimes \CD$.
\item For any finite category $\CX$, the functor
\[
- \circ \boxtimes \colon \fun_\bk^\rex(\CC,\CD;\CX) \to \fun_\bk^\rex(\CC \boxtimes \CD,\CX)
\]
is an equivalence.
\eit
\end{defn}

If $\CC \simeq \LMod_A(\vect_\bk)$ and $\CD \simeq \LMod_B(\vect_\bk)$ for some finite-dimensional $\bk$-algebra $A,B$, we have
\[
\CC \boxtimes_\bk \CD \simeq \LMod_{A \otimes_\bk B}(\vect_\bk) .
\]
In particular, if $\CC$ and $\CD$ are finite separable, their Deligne tensor product $\CC \boxtimes_\bk \CD$ is also finite separable. Therefore, Deligne tensor product defines a symmetric monoidal structure on the 2-category $2\vect_\bk$ of finite separable categories, $\bk$-linear functors and natural transformations. Another important property of Deligne tensor product is the following: for any $x,x' \in \CC$ and $y,y' \in \CD$, there is an isomorphism of $\bk$-vector spaces:
\[
(\CC \boxtimes \CD)(x \boxtimes y,x' \boxtimes y') \simeq \CC(x,x') \otimes_\bk \CD(y,y') .
\]
When $x = x'$ and $y = y'$, this is an isomorphism of $\bk$-algebras.

\begin{defn}
A \emph{multi-fusion category over $\bk$} (or a \emph{$\bk$-linear multi-fusion category}) is a rigid monoidal $\bk$-linear category that is finite separable. It is \emph{fusion} if the tensor unit $\one$ is a simple object.
\end{defn}

\begin{rem}
A multi-fusion category over $\bk$ is the same as a rigid algebra (see for example \cite{Dec23a}) in $2\vect_\bk$.
\end{rem}

\begin{expl}
Let $G$ be a finite group. The category $\rep_\bk(G)$ of finite-dimensional $G$-representations over $\bk$ is a (symmetric) fusion category if the characteristic of $\bk$ does not divide the order of $G$. The category $\vect_\bk[G]$ of finite-dimensional $G$-graded vector spaces over $\bk$ is a fusion category.
\end{expl}

Let $\CC$ be a fusion category over $\bk$. Define $\Omega \CC \coloneqq \CC(\one,\one)$. By the Eckmann-Hilton argument, $\Omega \CC$ is a field extension of $\bk$ of finite degree. For any $x,y \in \CC$, there are two $\Omega \CC$-actions on the hom space $\CC(x,y)$:
\[
\begin{aligned}
\lambda \triangleright_l f & \coloneqq \bigl( x \simeq \one \otimes x \xrightarrow{\lambda \otimes f} \one \otimes y \simeq y \bigr) , \\
\lambda \triangleright_r f & \coloneqq \bigl( x \simeq x \otimes \one \xrightarrow{f \otimes \lambda} y \otimes \one \simeq y \bigr) ,
\end{aligned} \quad \lambda \in \Omega \CC , \, f \in \CC(x,y) .
\]
If these two actions coincide, $\CC$ is a fusion category over $\Omega \CC$.

\begin{lem} \label{lem_braided_fusion_category_over_Omega}
Suppose $\CC$ admits a braiding. Then $\CC$ is a fusion category over $\Omega \CC$.
\end{lem}

\pf
The naturality of the braiding implies that the two $\Omega \CC$ actions on each hom space coincide.
\epf


\subsection{Modules over and algebras in multi-fusion categories}

Let $\CC$ be a multi-fusion category over $\bk$. A \emph{finite left $\CC$-module} is a finite $\bk$-linear category $\CM$ equipped with a left $\CC$-action such that the action functor $\odot \colon \CC \times \CM \to \CM$ is right exact and $\bk$-linear in each variable. Then for any $x \in \CM$, since $- \odot x \colon \CC \to \CM$ is right exact and both $\CC$ and $\CM$ are finite, it has a right adjoint $[x,-] \colon \CM \to \CC$ (called the \emph{internal hom} functor). When $x$ is a projective generator (the existence is provided by the finiteness of $\CM$), $[x,-]$ is faithful exact and thus monadic. The monad associated to the adjunction $(- \odot x) \dashv [x,-]$ is $[x,- \odot x]$, which is isomorphic to $- \odot [x,x]$ by the rigidity of $\CC$. So the monadicity theorem implies the following result.

\begin{thm}[\textnormal{\cite[Theorem 1]{Ost03}}] \label{thm_Ostrik_module_algebra}
Every finite left $\CC$-module $\CM$ is equivalent to the category $\RMod_A(\CC)$ of right $A$-modules in $\CC$ for some algebra $A \in \CC$.
\end{thm}

If $\CM \simeq \RMod_A(\CC)$ and $\CN \simeq \RMod_B(\CC)$ for some algebra $A,B \in \CC$, the category $\fun_\CC^\rex(\CM,\CN)$ of right exact left $\CC$-module functor $\CM \to \CN$ is equivalent to the category $\BMod_{A|B}(\CC)$ of $A$-$B$-bimodules in $\CC$. As a consequence, the 2-category of finite left $\CC$-modules, right exact left $\CC$-module functors and left $\CC$-module natural transformations is equivalent to the 2-category of algebras, bimodules and bimodule morphisms in $\CC$.

\begin{defn}
An algebra $A \in \CC$ is \emph{semisimple} if $\RMod_A(\CC)$ is semisimple. An algebra $A \in \CC$ is \emph{separable} if there exists an $A$-$A$-bimodule morphism $A \to A \otimes A$ that is a section of the multiplication of $A$.
\end{defn}

By the definition of a separable algebra, it is easy to prove the following result.

\begin{lem}
Let $A \in \CC$ be a separable algebra. Then for every semisimple left $\CC$-module $\CM$, the category $\LMod_A(\CM)$ is semisimple. Similarly, for every semisimple right $\CC$-module $\CN$, the category $\RMod_A(\CN)$ is semisimple. In particular, $A$ is semisimple.
\end{lem}

\begin{thm}[\textnormal{\cite[Theorem 3.5.4]{DSPS20} \cite[Proposition 4.6]{KZ17}}] \label{thm_separable_module}
Let $A \in \CC$ be an algebra and $\CM \coloneqq \RMod_A(\CC)$. Then the following conditions are equivalent:
\bnu[(1)]
\item $A$ is separable.
\item $A$ is Morita equivalent to a separable algebra, i.e., $\CM \simeq \RMod_B(\CC)$ for some separable algebra $B \in \CC$.
\item For every semisimple left $\CC$-module $\CN$, the category $\fun_\CC^\rex(\CM,\CN)$ is semisimple.
\item The category $\fun_\CC^\rex(\CM,\CM)$ is semisimple.
\item The category $\BMod_{A|A}(\CC)$ is semisimple.
\enu
\end{thm}

\begin{defn}
A finite left $\CC$-module is \emph{separable} if it is equivalent to $\RMod_A(\CC)$ for some separable algebra $A \in \CC$.
\end{defn}

\begin{rem} \label{rem_separable_module_separable_category}
Let $\CM$ be a separable finite left $\CC$-module. Then the underlying finite category of $\CM$ is separable \cite[Corollary 4.5]{KZ17}. The converse is not true. For example, if $G$ is a finite group and the characteristic of $\bk$ divides the order of $G$, then $\vect_\bk$ is not a separable left $\vect_\bk[G]$-module, because $\fun_{\vect_\bk[G]}(\vect_\bk,\vect_\bk) \simeq \rep_\bk(G)$ is not semisimple.
\end{rem}

Therefore, the 2-category $\LMod_\CC^\sep(2\vect_\bk)$ of separable left $\CC$-modules, left $\CC$-module functors and left $\CC$-module natural transformations is equivalent to the 2-category $\Sigma \CC$ of separable algebras, bimodules and bimodule morphisms in $\CC$. It is a \emph{compact semisimple 2-category} \cite[Theorem 1.3.4]{Dec23}.

\begin{rem} \label{rem_char_0_separable_semisimple}
Suppose the base field $\bk$ has characteristic zero. Then every finite semisimple left $\CC$-module is separable \cite[Corollary 3.6.9]{DSPS20}. Equivalently, every semisimple algebra in $\CC$ is separable \cite[Theorem 6.10]{KZ17}.
\end{rem}

\subsection{Morita theory of fusion categories}

Let $\CB,\CC,\CD$ be multi-fusion categories over $\bk$. Suppose $\CM$ is a finite right $\CC$-module, $\CN$ is a finite left $\CC$-module and $\CP$ is a finite category. A \emph{$\CC$-balanced functor} $\CM \times \CN \to \CP$ is a functor $F \colon \CM \times \CN \to \CP$ equipped with a natural isomorphism
\[
F(m \odot c,n) \simeq F(m,c \odot n) , \quad m \in \CM , \, n \in \CN , \, c \in \CC ,
\]
satisfying some coherence conditions. The category of right exact $\CC$-balanced functors $\CM \times \CN \to \CP$ is denoted by $\fun_\CC^\rex(\CM,\CN;\CP)$.

\begin{defn}
The \emph{relative tensor product} $\CM \boxtimes_\CC \CN$ is a finite category defined by the following universal property:
\bit
\item It is equipped with a right exact $\CC$-balanced functor $\boxtimes_\CC \colon \CM \times \CN \to \CM \boxtimes_\CC \CN$.
\item For any finite category $\CX$, the functor
\[
- \circ \boxtimes_\CC \colon \fun_\CC^\rex(\CM,\CN;\CX) \to \fun_\bk^\rex(\CM \boxtimes_\CC \CN,\CX)
\]
is an equivalence.
\eit
\end{defn}

If $\CM \simeq \LMod_A(\CC)$ and $\CN \simeq \RMod_B(\CC)$ for some algebra $A,B \in \CC$, we have
\[
\CM \boxtimes_\CC \CN \simeq \BMod_{A|B}(\CC) .
\]
When $\CC = \vect_\bk$, the relative tensor product $\boxtimes_{\vect_\bk}$ is the Deligne tensor product $\boxtimes_\bk$.

In addition, if $\CM$ is a finite $\CB$-$\CC$-bimodule and $\CN$ is a finite $\CC$-$\CD$-bimodule, the relative tensor product $\CM \boxtimes_\CC \CN$ is a finite $\CB$-$\CD$-bimodule.

\begin{thm}[\textnormal{\cite[Theorem 3.5.5]{DSPS20}}] \label{thm_separable_module_relative_tensor}
Suppose $\CM$ is a separable finite $\CB$-$\CC$-bimodule and $\CN$ is a separable finite $\CC$-$\CD$-bimodule. The relative tensor product $\CM \boxtimes_\CC \CN$ is a separable finite $\CB$-$\CD$-bimodule.
\end{thm}

\begin{rem}
Let $\CM$ be a separable finite left $\CC$-module. Then $\CM \simeq \CC \boxtimes_\CC \CM$ is a separable finite left $\vect_\bk$-module. In other words, $\CM \in 2\vect_\bk$ is a finite separable category (see Remark \ref{rem_separable_module_separable_category}).
\end{rem}

Hence, multi-fusion categories, separable bimodules, bimodule functors and bimodule natural transformations form a 3-category
. The equivalence of objects in this 3-category is Morita equivalence. Now we give a precise definition.

Let $\CM$ be a finite $\CC$-$\CD$-bimodule. The opposite category $\CM^\op$ admits two right $\CC$-actions:
\[
\begin{aligned}
x \odot^L c \coloneqq c^L \odot x , \\
x \odot^R c \coloneqq c^R \odot x ,
\end{aligned} \quad c \in \CC , \, x \in \CM .
\]
We denote these two right $\CC$-modules by $\CM^{\op|L}$ and $\CM^{\op|R}$. Similarly, $\CM^\op$ admits two left $\CD$-actions:
\[
\begin{aligned}
d \odot^L x \coloneqq x \odot d^L , \\
d \odot^R x \coloneqq x \odot d^R ,
\end{aligned} \quad d \in \CD , \, x \in \CM .
\]
We denote these two left $\CD$-modules by $\CM^{L|\op}$ and $\CM^{R|\op}$. If $\CM \simeq \RMod_A(\CC)$ for some algebra $A \in \CC$ and $\CM \simeq \LMod_B(\CD)$ for some algebra $B \in \CD$, then
\[
\CM^{\op|L} \simeq \LMod_A(\CC) , \quad \CM^{R|\op} \simeq \RMod_B(\CD) .
\]
For simplicity, we denote $\CM^{R|\op|L}$ by $\CM^\op$.

\begin{thm} \label{thm_invertible_bimodule}
For any finite $\CC$-$\CD$-bimodule $\CM$, the following conditions are equivalent:
\bnu[(1)]
\item $\CM \boxtimes_\CD \CM^\op$ is equivalent to $\CC$ as $\CC$-$\CC$-bimodules;
\item $\CM^\op \boxtimes_\CC \CM$ is equivalent to $\CD$ as $\CD$-$\CD$-bimodules;
\item the canonical monoidal functor $\CC \to \fun_{\CD^\rev}(\CM,\CM)$ is an equivalence;
\item the canonical monoidal functor $\CD^\rev \to \fun_\CC(\CM,\CM)$ is an equivalence;
\item there exists a finite $\CD$-$\CC$-bimodule $\CN$ such that $\CM \boxtimes_\CD \CN \simeq \CC$ as $\CC$-$\CC$-bimodules and $\CN \boxtimes_\CC \CM \simeq \CD$ as $\CD$-$\CD$-bimodules.
\enu
\end{thm}

The equivalence of the first four conditions is proved in \cite[Proposition 4.2]{ENO10}, and the equivalence between (3) and (4) is the double centralizer theorem \cite[Theorem 3.27]{EO04}. The last condition means that $\CM$ is an invertible 1-morphism in the 3-category of multi-fusion categories, finite bimodules, right exact bimodule functors and bimodule natural transformations. It is well-known that every 1-morphism in this 3-category is dualizable: the left dual of a finite $\CC$-$\CD$-bimodule $\CM$ is $\CM^{R|\op|R}$, and the right dual is $\CM^{L|\op|L}$. If $\CM$ is invertible, its inverse must be the dual. In this case, the different $\CD$-$\CC$-bimodule structures on $\CM^\op$ are equivalent.

\begin{defn}
A finite $\CC$-$\CD$-bimodule is \emph{invertible} if it satisfies one (and hence all) of the conditions in Theorem \ref{thm_invertible_bimodule}. We say $\CC$ and $\CD$ are \emph{Morita equivalent} if there exists an invertible $\CC$-$\CD$-bimodule.
\end{defn}

By Theorem \ref{thm_separable_module} (4) and Theorem \ref{thm_invertible_bimodule} (3) (4), every invertible $\CC$-$\CD$-bimodule must be separable. The following result is obvious.

\begin{thm} \label{thm_Morita_equivalent_defined_by_2-category}
Suppose $\CC$ and $\CD$ are Morita equivalent with invertible bimodules $\CM$ and $\CN$. The 2-categories $\LMod_\CC^\sep(2\vect_\bk)$ and $\LMod_\CD^\sep(2\vect_\bk)$ are equivalent via the 2-functors $\CM \boxtimes_\CD -$ and $\CN \boxtimes_\CC -$.
\end{thm}

Also by Theorem \ref{thm_invertible_bimodule}, if $\CC$ and $\CD$ are Morita equivalent, there exists a separable algebra $A \in \CC$ such that $\CD \simeq \BMod_{A|A}(\CC)$.

\subsection{Equivariantization}

\begin{defn}
Let $G$ be a group and $\CC$ be a category. A \emph{G-action} on $\CC$ is a monoidal functor $T \colon G \to \Aut(\CC)$, where $\Aut(\CC)$ is the 2-group of auto-equivalences of $\CC$ and natural isomorphisms. In other words, a $G$-action on $\CC$ consists of
\bit
\item auto-equivalences $T_g \colon \CC \to \CC$ for all $g \in G$,
\item and natural isomorphisms $\gamma_{g,h} \colon T_g \circ T_h \Rightarrow T_{gh}$ for all $g,h \in G$,
\eit
such that the following diagram commutes:
\[
\xymatrix{
T_g \circ T_h \circ T_k \ar@{=>}[r]^-{\gamma_{g,h} \circ 1} \ar@{=>}[d]_{1 \circ \gamma_{h,k}} & T_{gh} \circ T_k \ar@{=>}[d]^{\gamma_{gh,k}} \\
T_g \circ T_{hk} \ar@{=>}[r]^-{\gamma_{g,hk}} & T_{ghk}
}
\]
\end{defn}

\begin{defn}
Let $G$ be a group and $\CC$ be a category equipped with a $G$-action. The \emph{equivariantization} (or \emph{homotopy fixed points}) of $\CC$ is the category $\CC^G$ defined as follows:
\bit
\item The objects in $\CC^G$ are pairs $(x,u = \{u_g\}_{g \in G})$, where $x \in \CC$ and $u_g \colon g \odot x \to x$ is an isomorphism in $\CC$ such that for any $g,h \in G$ the following diagrams commute:
\be \label{diag_equivariantization_object}
\begin{array}{c}
\xymatrix{
T_g T_h(x) \ar[r]^{(\gamma_{g,h})_x} \ar[d]_{1 \odot u_h} & T_{gh}(x) \ar[d]^{u_{gh}} \\
T_g(x) \ar[r]^-{u_g} & x
}
\end{array}
\ee
\item A morphism $f \colon (x,u) \to (y,v)$ is a morphism $f \colon x \to y$ in $\CC$ such that for any $g \in G$ the following diagram commutes:
\be \label{diag_equivariantization_morphism}
\begin{array}{c}
\xymatrix{
T_g(x) \ar[r]^{1 \odot f} \ar[d]_{u_g} & T_g(y) \ar[d]^{v_g} \\
x \ar[r]^{f} & y
}
\end{array}
\ee
\item The composition and identity morphisms are induced by those of $\CC$.
\eit
Note that there is a forgetful functor $\CC^G \to \CC$ sending $(x,u)$ to $x$.
\end{defn}

\begin{rem}
The commutative diagram \eqref{diag_equivariantization_object} already implies that $u_g$ is an isomorphism for every $g \in G$.
\end{rem}

\begin{rem}
When $\CC$ is an $E_k$-monoidal category and the $G$-action on $\CC$ is given by a monoidal functor $G \to \Aut^{E_k}(\CC)$, then the equivariantization $\CC^\CG$ is also an $E_k$-monoidal category, and the forgetful functor $\CC^\CG \to \CC$ is an $E_k$-monoidal functor. Indeed, the tensor product of two objects $(x,u),(y,v) \in \CC^\CG$ is the object $x \otimes y \in \CC$ equipped with the isomorphisms $g \odot (x \otimes y) \simeq (g \odot x) \otimes (g \odot y) \xrightarrow{u_g \otimes v_g} x \otimes y$.
\end{rem}

Now we assume that $G$ is a finite group.

\begin{lem} \label{lem_equivariantization_monadic_comonadic}
Let $\CC$ be an additive category equipped with a $G$-action. Then the forgetful functor $\CC^G \to \CC$ is both monadic and comonadic.
\end{lem}

\pf
First we remark that for any $G$-action $T \colon G \to \Aut(\CC)$, the functor $T_g$ for every $g \in G$ is an equivalence and hence preserves direct sums. So the $G$-action is automatically an `additive action'.

The monoidal structure on $T$ induces both a monad and a comonad structures on the endofunctor
\[
\bigoplus_{g \in G} T_g \colon \CC \to \CC .
\]
Indeed, this is a Frobenius monad. By definition, a $(\bigoplus_g T_g)$-module in $\CC$ is an object $x \in \CC$ equipped with a morphism $u \colon \bigoplus_g T_g(x) \to x$ satisfying the associativity condition. The morphism $u$ is equivalent to a collection of morphisms $\{u_g \colon T_g(x) \to x\}_{g \in G}$, and the associativity is equivalent to the commutative diagram \eqref{diag_equivariantization_object}. Similarly, a $(\bigoplus_g T_g)$-module morphism $(x,u) \to (y,v)$ is a morphism $f \colon x \to y$ such that the diagram \eqref{diag_equivariantization_morphism} commutes. Therefore, the category of $(\bigoplus_g T_g)$-modules (Eilenberg-Moore category of $\tilde T$) is equivalent to $\CC^G$. One can similarly show that $\CC^G$ is also equivalent to the category of $(\bigoplus_g T_g)$-comodules.
\epf

\begin{expl}
For the trivial $G$-action on $\vect_\bk$, the associated Frobenius monad is $\bk[G] \otimes_\bk -$, where $\bk[G]$ is equipped with the following Frobenius algebra structure:
\bit
\item The algebra structure on $\bk[G]$ is the usual group algebra structure.
\item The coalgebra structure on $\bk[G]$ is defined by
\[
\Delta(g) \coloneqq \sum_{h \in G} gh \otimes h^{-1} , \quad \varepsilon(g) = \delta_{g,e} \cdot 1 .
\]
\eit
Both the category of $\bk[G]$-modules and $\bk[G]$-comodules are equivalent to $\rep_\bk(G) \simeq \vect_\bk^G$.
\end{expl}

Let $\CC$ be a finite $\bk$-linear category equipped with a $\bk$-linear $G$-action. Clearly this action is equivalent to a left $\vect_\bk[G]$-module structure on $\CC$.

\begin{lem} \label{lem_equivariantization_module_group_algebra}
The equivariantization $\CC^G$ is equivalent to $\LMod_{\bk[G]}(\CC)$ as $\bk$-linear categories. In particular, $\CC^G$ is also finite.
\end{lem}

\pf
By Lemma \ref{lem_equivariantization_monadic_comonadic}, $\CC^G$ is equivalent to the category of modules over the monad $\bigoplus_g T_g$. The action functor $T_g$ is isomorphic to $\bk_g \odot - \colon \colon \CC \to \CC$, where $\{\bk_g\}_{g \in G}$ are simple objects of $\vect_\bk[G]$. Then one can check that
\[
\bigoplus_{g \in G} T_g \simeq \bigoplus_{g \in G} \bk_g \odot - \simeq \bk[G] \odot -
\]
as monads, where the monad structure on the right hand side is given by the algebra structure on the group algebra $\bk[G]$. This proves that $\CC^G \simeq \LMod_{\bk[G]}(\CC)$.
\epf

\section{Symmetric fusion categories over algebraically closed fields of characteristic zero} \label{sec_Deligne_symmetric_fusion}

In this section, we review the classification of symmetric fusion categories over an algebraically closed field of characteristic zero.

\subsection{Super groups and super representations}

In this subsection, we recall the construction of two families of symmetric fusion categories over any field $\bk$.

The first family is given by group representations. For any finite group $G$ such that the characteristic of $\bk$ does not divide the order of $G$, the category $\rep_\bk(G) \simeq \fun(\rmB G,\vect_\bk)$ of finite-dimensional $G$-representations over $\bk$ is a fusion category. It admits an obvious symmetric braiding inherited from $\vect_\bk$. We say that a symmetric fusion category is \emph{bosonic} if it is equivalent to $\rep_\bk(G)$ for some finite group $G$.

The second family is given by the super version of group representations. A \emph{super $\bk$-vector space} is a $\bk$-vector space $V$ equipped with a $\Zb_2$-grading: $V = V_0 \oplus V_1$. The category of finite-dimensional super $\bk$-vector spaces and grading-preserving $\bk$-linear maps is denoted by $\svect_\bk$. It is a symmetric fusion category with the following braiding structure:
\begin{align*}
V \otimes W & \to W \otimes V \\
v \otimes w & \mapsto (-1)^{mn} \cdot w \otimes v ,
\end{align*}
where $v \in V_m$ and $w \in W_n$. For any $V \in \svect$, the \emph{parity morphism} is defined as follows:
\[
\Pi_V \coloneqq \biggl( V = V_0 \oplus V_1 \xrightarrow{1 \oplus (-1)} V_0 \oplus V_1 = V \biggr) .
\]
Then $\Pi \colon 1_{\svect_\bk} \to 1_{\svect_\bk}$ is a monoidal natural automorphism. It is easy to check that the 2-group $\Aut^{E_3}_\bk(\svect_\bk)$ of $\bk$-linear symmetric monoidal autoequivalences of $\svect_\bk$ and monoidal natural isomorphisms is equivalent to $\mathrm B \Zb_2$:
\bit
\item Every $\bk$-linear symmetric monoidal autoequivalence of $\svect_\bk$ is monoidally isomorphic to the identity functor.
\item There is only one nontrivial monoidal natural automorphism of the identity functor, that is, the parity natural isomorphism $\Pi$.
\eit

\begin{defn}
A \emph{super groupoid} is a groupoid equipped with a $\rmB \Zb_2$-action. A \emph{super group} is a super groupoid with only one object.
\end{defn}

Explicitly, a super group is a group $G$ equipped with a distinguished central element $z \in G$ such that $z^2 = e$.

\begin{expl}
Let $V \in \svect_\bk$. The group of grading-preserving $\bk$-linear automorphisms of $V$ is denoted by $\mathrm{GL}_\bk(V)_0$. Then $(\mathrm{GL}_\bk(V)_0,\Pi_V)$ is a super group.
\end{expl}

\begin{defn}
Let $(G,z)$ and $(H,w)$ be super groups. A \emph{super group homomorphism} $f \colon (G,z) \to (H,w)$ is a group homomorphism $f \colon G \to H$ such that $f(z) = w$.
\end{defn}

Let $(G,z)$ be a super group.

\begin{defn}
A \emph{super representation (over $\bk$)} of $(G,z)$ is a super $\bk$-vector space $V = V_0 \oplus V_1$ equipped with a super group homomorphism $\rho \colon (G,z) \to (\mathrm{GL}_\bk(V)_0,\Pi_V)$. In other words, $(V,\rho)$ is a representation of $G$ such that $\rho(z) = \Pi_V$. A \emph{homomorphism} of super representations of $(G,z)$ is a homomorphism of the underlying representations of $G$.
The category of finite-dimensional super representations of $(G,z)$ is denoted by $\srep_\bk(G,z)$.
\end{defn}

The following lemma is almost obvious from the definition.

\begin{lem}
The category $\srep_\bk(G,z)$ is equivalent to the category $\fun_{\rmB \Zb_2}(\rmB G,\svect_\bk)$ of $\rmB \Zb_2$-equivariant functors from $\rmB G$ to $\svect_\bk$.
\end{lem}

\begin{lem}
The functor $\srep_\bk(G,z) \to \rep_\bk(G)$ defined by forgetting the grading is an equivalence.
\end{lem}

\pf
By definition this functor is fully faithful. Now we show that it is essentially surjective. Let $(V,\rho) \in \rep_\bk(G)$. Then $\rho(z)$ is a $\bk$-linear automorphism of $V$ satisfying $\rho(z)^2 = 1_V$. So $V$ can be decomposed into $V = V_0 \oplus V_1$, where
\[
V_j = \{v \in V \mid \rho(z)(v) = (-1)^j \cdot v\}
\]
are the eigenspaces of $\rho(z)$. Clearly $V = V_0 \oplus V_1$ is a super representation of $(G,z)$.
\epf

Therefore, $\srep_\bk(G,z)$ is a fusion category if $G$ is finite and the characteristic of $\bk$ does not divide the order of $G$. It admits a symmetric braiding inherited from $\svect_\bk$. If $z = e \in G$ is trivial, we have the equivalence of symmetric fusion categories
\[
\srep_\bk(G,z) \simeq \fun_{\rmB \Zb_2}(\rmB G,\svect_\bk) \simeq \fun(\rmB G,\svect_\bk^{\rmB \Zb_2}) \simeq \fun(\rmB G,\vect_\bk) \simeq \rep_\bk(G) .
\]
We say that a symmetric fusion category is \emph{fermionic} if it is equivalent to $\srep_\bk(G,z)$ for some finite super group $(G,z)$ where $z$ is nontrivial.

\subsection{Deligne's classification results}

In this subsection we review Deligne's results on the classification of symmetric fusion categories over an algebraically closed field of characteristic zero.

Let $\bk$ be a field.

\begin{defn}
A \emph{fiber functor} is a faithful exact $\bk$-linear symmetric monoidal functor. A symmetric fusion category $\CC$ over $\bk$ is \emph{Tannakian} if it admits a fiber functor $\CC \to \vect_\bk$.
\end{defn}

\begin{expl}
Let $(G,z)$ be a finite super group. The forgetful functor $\srep_\bk(G,z) \to \svect_\bk$ is a fiber functor.
\end{expl}

In \cite{DM82}, Deligne and Milne studied general Tannakian tensor categories (which are not necessarily finite nor semisimple) and fiber functors (which are required to be exact). In the case we are interested in, Tannakian symmetric fusion categories are classified as follows.

\begin{thm}
Suppose $\bk$ is algebraically closed. Let $\CC$ be a Tannakian symmetric fusion category and $\omega \colon \CC \to \vect_\bk$ be a fiber functor.
\bnu[(1)]
\item The group $G \coloneqq \Aut^{E_3}(\omega)$ of monoidal automorphisms of $\omega$ is a finite group.
\item The obvious lifting $\CC \to \rep_\bk(G)$ of $\omega$ is an equivalence.
\enu
\end{thm}

In other words, a symmetric fusion category over an algebraically closed field is bosonic if and only if it is Tannakian. This result is essentially a special case of \cite[Theorem 2.11]{DM82}. In general, not all symmetric fusion categories admit fiber functors to $\vect_\bk$. In \cite[Section 8.19]{Del90}, Deligne studied fiber functors to $\svect_\bk$.

\begin{thm} \label{thm_Deligne_symmetric_fusion_super_lifting}
Suppose $\bk$ is algebraically closed. Let $\CC$ be a symmetric fusion category and $\omega \colon \CC \to \svect_\bk$ be a fiber functor.
\bnu[(1)]
\item The group $G \coloneqq \Aut^{E_3}(\omega)$ is a finite group, and $z \coloneqq \Pi \circ 1_\omega \in G$ is a central element satisfying $z^2 = 1$.
\item The obvious lifting $\CC \to \srep_\bk(G,z)$ of $\omega$ is an equivalence.
\enu
\end{thm}

In other words, a symmetric fusion category over an algebraically closed field is either bosonic or fermionic if and only if it admits a fiber functor to $\svect_\bk$. In \cite{Del02}, Deligne studied the existence of fiber functors to $\svect_\bk$ over algebraically closed fields of characteristic zero.

\begin{thm}[\textnormal{\cite[Corollary 0.8]{Del02}}] \label{thm_Deligne_classificaion_symmetric_fusion_category}
Suppose $\bk$ is algebraically closed and of characteristic zero. Then every symmetric fusion category $\CC$ admits a fiber functor to $\svect_\bk$. Therefore, it is equivalent to $\srep_\bk(G,z)$ for some finite super group $(G,z)$.
\end{thm}

\begin{rem}
Suppose $\bk$ is algebraically closed and the characteristic of $\bk$ is $p > 0$. Then every symmetric fusion category admits a fiber functor to the Verlinde category $\mathrm{Ver}_p$ \cite[Theorem 1.5]{Ost20}.
\end{rem}

Now let $\bk$ be an algebraically closed field of characteristic zero and $\CC$ be a symmetric fusion category over $\bk$. We use $\SF_\CC$ to denote the groupoid of fiber functors $\CC \to \svect_\bk$ and monoidal natural isomorphisms. It admits an action of $\Aut^{E_3}_\bk(\svect_\bk) \simeq \rmB \Zb_2$.

\begin{thm} \label{thm_uniqueness_fiber_functor}
All fiber functors $\CC \to \svect_\bk$ are monoidally isomorphic to each other. In other words, $\SF_\CC$ is connected.
\end{thm}

This result is \cite[Theorem 9.9.26 (ii)]{EGNO15} (see also \cite[Section 9.11]{EGNO15} or \cite[Section 6.4]{Cou20}). When $\CC = \rep_\bk(G)$ or $\srep_\bk(G,z)$, the groupoid $\SF_\CC$ is also denoted by $\SF_G$ or $\SF_{(G,z)}$, respectively. By the classical Tannaka-Krein duality for finite groups, the functor $\rmB G \to \SF_{(G,z)}$ sending the unique object to the forgetful functor $\srep_\bk(G,z) \to \svect_\bk$ is a $\rmB \Zb_2$-equivariant equivalence.

\subsection{Autoequivalences of symmetric fusion categories}

In this subsection, $\bk$ is an algebraically closed field of characteristic zero.


Let $G$ be a finite group. All autoequivalences of the groupoid $\rmB G$ and natural isomorphisms form a 2-group $\Aut(\rmB G)$.
\bit
\item The objects of $\Aut(\rmB G)$ are one-to-one corresponding to the group automorphisms of $G$. The tensor product of objects is given by the composition of automorphisms.
\item Given two group automorphisms $\phi$ and $\psi$, viewed as autoequivalences of $\mathrm B G$, the natural isomorphisms $\phi \Rightarrow \psi$ are given by the elements $g \in G$ satisfying $\psi(x) g = g \phi(x)$ for all $x \in G$. The composition of morphisms is given by the multiplication of $G$.
\eit
All $\bk$-linear symmetric autoequivalences of $\rep_\bk(G)$ and monoidal natural isomorphisms also form a 2-group $\Aut_\bk^{E_3}(\rep_\bk(G))$.

We identify $\rep_\bk(G)$ with the functor category $\fun(\rmB G,\vect_\bk)$. Then for any $\phi \in \Aut(\rmB G)$, there is an equivalence $\phi^*$ of $\rep_\bk(G)$ defined as follows:
\begin{align*}
\phi^* \colon \rep_\bk(G) & \to \rep_\bk(G) \\
(V,\rho) & \mapsto (V,\rho \circ \phi) .
\end{align*}
Moreover, since the symmetric monoidal structure of $\rep_\bk(G)$ is inherited from $\vect$, the autoequivalence $\phi^*$ is also symmetric monoidal. Clearly $\phi^* \circ \psi^* = (\psi \circ \phi)^*$ for any $\phi,\psi \in \Aut(\rmB G)$. So $\phi \mapsto (\phi^{-1})^*$ defines a monoidal functor
\[
\Phi \colon \Aut(\mathrm B G) \to \Aut_\bk^{E_3}(\rep_\bk(G)) .
\]
To be more precise, this functor $\Phi$ sends a morphism $g \colon \phi \to \psi$ in $\Aut(\mathrm B G)$ to the monoidal natural isomorphism $\alpha^g \colon (\phi^{-1})^* \Rightarrow (\psi^{-1})^*$ defined by
\[
\alpha^g_{(V,\rho)} \coloneqq \rho(\psi^{-1}(g))^{-1} \colon (V,\rho \circ \phi^{-1}) \to (V,\rho \circ \psi^{-1}) .
\]

\begin{prop} \label{prop_auto_rep_G_linear}
The monoidal functor $\Phi \colon \Aut(\rmB G) \to \Aut_\bk^{E_3}(\rep_\bk(G))$ is an equivalence.
\end{prop}

\pf
As a monoidal functor between 2-groups, $\Phi$ is fully faithful if and only if its map on the second homotopy group is an isomorphism. It is easy to see that
\[
\pi_2 (\Aut(\mathrm B G)) = \Aut(1_{\mathrm B G}) \simeq Z(G)
\]
and
\[
\pi_2 (\Aut_{E_3}(\rep_\bk(G))) = \Aut_{E_3}(1_{\rep_\bk(G)}) \simeq Z(G) ,
\]
where the last isomorphism maps $z \in Z(G)$ to the monoidal natural isomorphism $\alpha^z$ defined by $\alpha^z_{(V,\rho)} = \rho(z)$. This shows that $\Phi$ is fully faithful.

Then we show that $\Phi$ is essentially surjective. Let $\SF_G$ be the groupoid of fiber functors $\rep_\bk(G) \to \vect_\bk$. There is a functor $\ev \colon \rep_\bk(G) \to \fun(\SF_G,\vect_\bk)$ defined by $(V,\rho) \mapsto \ev_{(V,\rho)}$, where $\ev_{(V,\rho)}$ is defined by
\begin{align*}
\ev_{(V,\rho)} \colon \SF_G & \to \vect_\bk \\
U & \mapsto U(V,\rho) .
\end{align*}
One can check that the following diagram commutes:
\begin{equation} \label{diag_rep_fiber_functor}
\begin{array}{c}
\xymatrix{
 & \rep_\bk(G) \ar[dl]_{\ev} \ar@{<->}[dr]^{\simeq} \\
\fun(\SF_G,\vect_\bk) \ar[rr]^-{- \circ \iota} & & \fun(\rmB G,\vect_\bk)
}
\end{array}
\end{equation}
where $\iota \colon \rmB G \to \SF_G$ is the equivalence sending the unique object to the forgetful functor $\rep_\bk(G) \to \vect_\bk$. Therefore, $\ev \colon \rep_\bk(G) \to \fun(\SF_G,\vect_\bk)$ is also an equivalence.

Given a symmetric monoidal autoequivalence $F \in \Aut_\bk^{E_3}(\rep_\bk(G))$, there is an autoequivalence
\begin{align*}
F^* \colon \SF_G & \to \SF_G \\
U & \mapsto U \circ F^ .
\end{align*}
Then we define $\tilde F \coloneqq \iota^{-1} \circ (F^{-1})^* \circ \iota \in \Aut(\rmB G)$. By definition, the symmetric monoidal equivalence $\Phi(\tilde F)$ renders the following diagram commutative (in the 2-category of symmetric monoidal categories):
\[
\xymatrix@C=4em{
\rep_\bk(G) \ar@{-->}[r]^-{\Phi(\tilde F)} \ar@{<->}[d]_{\simeq} & \rep_\bk(G) \ar@{<->}[d]^{\simeq} \\
\fun(\rmB G,\vect_\bk) \ar[r]^-{- \circ \tilde{F}^{-1}} & \fun(\rmB G,\vect_\bk)
}
\]
Let us consider the following diagram (in the 2-category of symmetric monoidal categories):
\[
\xymatrix@C=4em{
\rep_\bk(G) \ar[r]^-{F} \ar@{<->}[d]_{\simeq} \ar@/_5em/[dd]_{\ev} & \rep_\bk(G) \ar@{<->}[d]^{\simeq} \ar@/^5em/[dd]^{\ev} \\
\fun(\rmB G,\vect_\bk) \ar[r]^-{- \circ \tilde{F}^{-1}} & \fun(\rmB G,\vect_\bk) \\
\fun(\SF_G,\vect_\bk) \ar[r]^-{- \circ F^*} \ar[u]^{- \circ \iota} & \fun(\SF_G,\vect_\bk) \ar[u]_{- \circ \iota}
}
\]
The left and right triangles are the commutative diagram \eqref{diag_rep_fiber_functor}. The commutativity of the outer diagram is obvious. The lower rectangle commutes by definition of $\tilde F$. So the upper rectangle also commutes, up to monoidal natural isomorphisms. Hence $\Phi(\tilde F)$ is monoidally isomorphic to $F$.
\epf


Now let $(G,z)$ be a finite super group. There is a 2-group $\Aut_{\rmB \Zb_2}(\rmB G)$ of $\rmB \Zb_2$-equivariant autoequivalences of $\rmB G$ and natural isomorphisms.
\bit
\item The objects of $\Aut_{\mathrm B \Zb_2}(\mathrm B G)$ are one-to-one corresponding to the group automorphisms $\phi$ of $G$ such that $\phi(z) = z$.
\item Every natural isomorphism is $\rmB \Zb_2$-equivariant. So the morphisms of $\Aut_{\mathrm B \Zb_2}(\mathrm B G)$ are the same as those of $\Aut(\mathrm B G)$. In other words, the forgetful functor $\Aut_{\mathrm B \Zb_2}(\mathrm B G) \to \Aut(\mathrm B G)$ is fully faithful.
\eit
Similar to the bosonic case, by identifying $\srep_\bk(G,z)$ with $\fun_{\rmB \Zb_2}(\rmB G,\svect_\bk)$, we obtain a monoidal functor $\Aut_{\rmB \Zb_2}(\rmB G) \to \Aut_\bk^{E_3}(\srep_\bk(G,z))$.

\begin{prop} \label{prop_auto_srep_G_z_linear}
The monoidal functor $\Aut_{\rmB \Zb_2}(\rmB G) \to \Aut_\bk^{E_3}(\srep_\bk(G,z))$ is an equivalence.
\end{prop}

\pf
The proof is parallel to Proposition \ref{prop_auto_rep_G_linear}. Let $\SF_{(G,z)}$ be the groupoid of fiber functors $\srep_\bk(G,z) \to \svect_\bk$. Then we also have a commutative diagram
\[
\xymatrix{
 & \srep_\bk(G,z) \ar[dl]_{\ev} \ar@{<->}[dr]^{\simeq} \\
\fun_{\rmB \Zb_2}(\SF_{(G,z)},\svect_\bk) \ar[rr]^-{- \circ \iota} & & \fun_{\rmB \Zb_2}(\rmB G,\svect_\bk)
}
\]
where $\iota \colon \rm B G \to \SF_{(G,z)}$ is the $\rmB \Zb_2$-equivariant equivalence sending the unique object to the forgetful functor $\srep_\bk(G,z) \to \svect_\bk$. The remainder of the proof proceeds similarly to the proof of Proposition \ref{prop_auto_rep_G_linear}, so we omit the details.
\epf

We denote the 2-groupoid of $\bk$-linear symmetric multi-fusion categories, $\bk$-linear symmetric monoidal equivalences and monoidal natural isomorphisms by $\Alg_{E_3}^\rig(2\vect_\bk)^\times$, and denote the 2-groupoid of finite super groupoids, $\rmB \Zb_2$-equivariant equivalences and natural isomorphisms by $\LMod_{\rmB \Zb_2}(\mathrm{fGrpd})^\times$.

\begin{thm} \label{thm_symmetric_fusion_category_C_super_groupoid}
The functors
\begin{align*}
(\Alg_{E_3}^\rig(2\vect_\bk)^\times)^\op & \to \LMod_{\rmB \Zb_2}(\mathrm{fGrpd})^\times & (\LMod_{\rmB \Zb_2}(\mathrm{fGrpd})^\times)^\op & \to \Alg_{E_3}^\rig(2\vect_\bk)^\times \\
\CC & \mapsto \SF_\CC & X & \mapsto \fun_{\rmB \Zb_2}(X,\svect_\bk)
\end{align*}
define an opposite equivalence between the 2-groupoids $\Alg_{E_3}^\rig(2\vect_\bk)^\times$ and $\LMod_{\rmB \Zb_2}(\mathrm{fGrpd})^\times$.
\end{thm}

\pf
Note that every symmetric multi-fusion category must be a direct sum of symmetric fusion categories. Then by Theorem \ref{thm_Deligne_classificaion_symmetric_fusion_category} the 2-functor $X \mapsto \fun_{\rmB \Zb_2}(X,\svect_\bk)$ is essentially surjective, and by Proposition \ref{prop_auto_srep_G_z_linear} it is fully faithful. Also by Theorem \ref{thm_Deligne_symmetric_fusion_super_lifting}, the quasi-inverse is given by $\CC \mapsto \SF_\CC$.
\epf

\begin{rem}
Theorem \ref{thm_symmetric_fusion_category_C_super_groupoid} also appeared in \cite[Theorem 2.79]{DHJ+24}.
\end{rem}

\section{Complexification and real forms of categories} \label{sec_Galois_descent}

In this section, we discuss the complexification and real forms of finite categories.

\subsection{Complexification}

Let $\CC$ be a finite $\Rb$-linear category.

\begin{defn}
A \emph{complexification} of $\CC$ is a finite $\Cb$-linear category $\CC_\Cb$ equipped with a right exact $\Rb$-linear functor $\iota \colon \CC \to \CC_\Cb$ satisfying the following universal property: for any $\Cb$-linear finite category $\CD$, the functor
\[
- \circ \iota \colon \fun_\Cb^{\mathrm{rex}}(\CC_\Cb,\CD) \to \fun_\Rb^{\mathrm{rex}}(\CC,\CD)
\]
is an equivalence, where $\fun_\Rb^{\mathrm{rex}}$ and $\fun_\Cb^{\mathrm{rex}}$ denote the category of $\Rb$-linear right exact functors and the category of $\Cb$-linear right exact functors, respectively.
\end{defn}

\begin{rem}
When $\CC$ is $E_k$-monoidal such that the tensor product is right exact and $\Rb$-linear in each variable, the complexification $\CC_\Cb$ is also $E_k$-monoidal and the functor $\CC \to \CC_\Cb$ is an $E_k$-monoidal functor.
\end{rem}

\begin{prop} \label{prop_complexification_algebra}
Let $A$ be a finite-dimensional $\Rb$-algebra. Denote $A_\Cb \coloneqq \Cb \otimes_\Rb A$. Then $\LMod_{A_\Cb}(\vect_\Cb)$ equipped with the functor $- \otimes_\Rb \Cb \colon \LMod_A(\vect_\Rb) \to \LMod_{A_\Cb}(\vect_\Cb)$ is a complexification of $\LMod_A(\vect_\Rb)$.
\end{prop}

\pf
Let $B$ be a finite-dimensional $\Cb$-algebra and $\CD \coloneqq \LMod_B(\vect_\Cb) \simeq \LMod_B(\vect_\Rb)$. By the Eilenberg-Watts theorem, there are equivalences
\[
\fun_\Cb^{\mathrm{rex}}(\LMod_{A_\Cb}(\vect_\Cb),\CD) \simeq \BMod_{B|A_\Cb}(\vect_\Cb) , \quad \fun_\Rb^{\mathrm{rex}}(\LMod_A(\vect_\Rb),\CD) \simeq \BMod_{B|A}(\vect_\Rb) .
\]
Note that the functor $- \otimes_\Rb \Cb \colon \LMod_A(\vect_\Rb) \to \LMod_{A_\Cb}(\vect_\Cb)$ is isomorphic to the functor
\[
A_\Cb \otimes_A - \colon \LMod_A(\vect_\Rb) \to \LMod_{A_\Cb}(\vect_\Rb) \simeq \LMod_{A_\Cb}(\vect_\Cb) .
\]
Under the above equivalences, the functor $\fun_\Cb^{\mathrm{rex}}(\LMod_{A_\Cb}(\vect_\Cb),\CD) \to \fun_\Rb^{\mathrm{rex}}(\LMod_A(\vect_\Rb),\CD)$ can be identified with the functor
\[
\BMod_{B|A_\Cb}(\vect_\Cb) \simeq \BMod_{B|A_\Cb}(\vect_\Rb) \xrightarrow{- \otimes_{A_\Cb} A_\Cb \otimes_A A} \BMod_{B|A}(\vect_\Rb) ,
\]
which is simply the forgetful functor. To show that this functor is an equivalence, consider the following commutative diagram of forgetful functors:
\[
\xymatrix{
\BMod_{B|A_\Cb}(\vect_\Cb) \ar[r] \ar[d] & \BMod_{B|A}(\vect_\Rb) \ar[d] \\
\LMod_B(\vect_\Cb) \ar[r]^-{\simeq} & \LMod_B(\vect_\Rb)
}
\]
Both two vertical arrows are monadic functors. The monad induced by the right vertical arrow is $(- \otimes_\Rb A)$, and the monad induced by the left vertical arrow is $(- \otimes_\Cb A_\Cb) \simeq (- \otimes_\Cb (\Cb \otimes_\Rb A)) \simeq (- \otimes_\Rb A)$. By the monadicity theorem, the above horizontal arrow is an equivalence.
\epf

\begin{cor}
The following categories are complexifications of $\CC$:
\bnu[(1)]
\item the Deligne tensor product $\vect_\Cb \boxtimes_\Rb \CC$ equipped with the functor $\Cb \boxtimes_\Rb - \colon \CC \to \vect_\Cb \boxtimes_\Rb \CC$;
\item the category $\LMod_\Cb(\CC)$ of left $\Cb$-modules in $\CC$ equipped with the free module functor $\Cb \odot - \colon \CC \to \LMod_\Cb(\CC)$, where $\CC$ is viewed as a left $\vect_\Rb$-module and $\Cb$ is viewed as an algebra in $\vect_\Rb$.
\enu
\end{cor}

In particular, if $\CC$ is finite semisimple, the complexification $\CC_\Cb$ is also finite semisimple.

\subsection{Semi-linear \texorpdfstring{$\Zb_2$}{Z2}-actions and Galois descent}

Recall that the complex conjugate of a complex number $\lambda = a + b \mathrm i \in \Cb$ is $\bar \lambda \coloneqq a - b \mathrm i \in \Cb$.

\begin{defn}
Let $V$ be a $\Cb$-vector space. Its \emph{complex conjugate} $\overline V$ is the $\Cb$-vector space defined by the same underlying abelian group of $V$ but equipped with the scalar multiplication
\begin{align*}
\Cb \times \overline V & \to \overline V \\
(\lambda , v) & \mapsto \bar \lambda \cdot v ,
\end{align*}
where $\bar \lambda \cdot v$ is the scalar multiplication in $V$.
\end{defn}

In other words, $\overline V$ is the pullback of $V$ along the ring homomorphism $\overline{(\cdot)} \colon \Cb \to \Cb$.

\begin{expl}
Let $V,W$ be $\Cb$-vector spaces. Then $\overline{\Hom_\Cb(V,W)}$ is isomorphic to $\Hom_\Cb(\overline V,\overline W)$ as $\Cb$-vector spaces via the identity map.
\end{expl}

It is clear that $\overline V \otimes \overline W = \overline{V \otimes W}$. Thus $V \mapsto \overline{V}$ defines a symmetric monoidal equivalence $\overline{( \cdot )} \colon \vect_\Cb \to \vect_\Cb$, which acts `trivially' on morphisms.


\begin{defn}
Let $V,W$ be $\Cb$-vector spaces. An \emph{anti-linear map} $f \colon V \to W$ is a homomorphism of abelian groups that satisfies
\[
f(\lambda \cdot v) = \bar \lambda \cdot f(v)
\]
for all $\lambda \in \Cb$ and $v \in V$.
\end{defn}

Equivalently, an anti-linear map $f \colon V \to W$ is a $\Cb$-linear map $f \colon \overline V \to W$.

\begin{defn}
Let $\CC$ be a $\Cb$-linear category. Its \emph{complex conjugate} $\overline \CC$ is defined by $\ob(\overline \CC) \coloneqq \ob(\CC)$ and $\Hom_{\overline \CC}(x,y) \coloneqq \overline{\Hom_\CC(x,y)}$.
\end{defn}

Similarly, $\CC \mapsto \overline{\CC}$ defines a symmetric monoidal equivalence $2\vect_\Cb \to 2\vect_\Cb$.

\begin{defn}
Let $\CC,\CD$ be $\Cb$-linear categories. A functor $F \colon \CC \to \CD$ is \emph{anti-linear} if the map $F_{x,y} \colon \Hom_\CC(x,y) \to \Hom_\CD(F(x),F(y))$ is anti-linear for every pair of objects $x,y \in \CC$.
\end{defn}

Equivalently, an anti-linear functor $F \colon \CC \to \CD$ is a $\Cb$-linear functor $F \colon \overline \CC \to \CD$. For example, the complex conjugation functor $\overline{( \cdot )} \colon \vect_\Cb \to \vect_\Cb$ is anti-linear.

\medskip
Let $\CD$ be a finite $\Cb$-linear category.

\begin{defn}
A \emph{semi-linear $\Zb_2$-action} is a $\Zb_2$-action on $\CD$ such that the generator of $\Zb_2$ acts as an anti-linear autoequivalence.
\end{defn}

Explicitly, a semi-linear $\Zb_2$-action on $\CD$ consists of
\bit
\item an anti-linear autoequivalence $J \colon \CD \to \CD$,
\item and a natural isomorphism $\phi \colon J \circ J \simeq 1_\CD$,
\eit
such that $J(\phi_x) = \phi_{J(x)} \in \CD(JJJ(x),J(x))$ for all $x \in \CD$.

\begin{expl} \label{expl_semi-linear_action_Vec}
The complex conjugation functor $\overline{(\cdot)} \colon \vect_\Cb \to \vect_\Cb$ with $\phi = \pm 1$ defines two different semi-linear $\Zb_2$-actions on $\vect_\Cb$. Indeed, the axioms for the $\Zb_2$-action imply that $\phi \in Z^2(\Zb_2;\Cb^\times_T)$ is a 2-cocycle, where $\Cb^\times_T$ is the group $\Cb^\times$ equipped with the $\Zb_2$-action given by complex conjugation. These two different actions correspond to the cohomology classes in $H^2(\Zb_2;\Cb^\times_T) \simeq \Zb_2$. For $\phi = 1$, the action is symmetric monoidal. For $\phi = -1$, the action is not monoidal.
\end{expl}

Now suppose $\CD$ is equipped with a semi-linear $\Zb_2$-action $(J,\phi)$ as above. We explicitly write down the defining data of the equivariantization $\CD^{\Zb_2}$.
\bit
\item The objects are pairs $(x,u)$, where $x \in \CD$ and $u \colon J(x) \simeq x$ is an isomorphism in $\CD$ such that the following diagram commutes:
\[
\xymatrix{
JJ(x) \ar[r]^-{J(u)} \ar[dr]_{\phi_x} & J(x) \ar[d]^{u} \\
 & x
}
\]
\item The morphisms $(x,u) \to (y,v)$ are morphisms $f \colon x \to y$ in $\CD$ such that the following diagram commutes:
\[
\xymatrix{
J(x) \ar[r]^-{J(f)} \ar[d]_{u} & J(y) \ar[d]^{v} \\
x \ar[r]^-{f} & y
}
\]
\eit

\begin{expl} \label{expl_semi-linear_equivariantization_Vec}
Let us compute the equivariantization $\vect_\Cb^{\Zb_2}$, where $\vect_\Cb$ is equipped with the semi-linear $\Zb_2$-action given in Example \ref{expl_semi-linear_action_Vec}.

First we consider the action with $\phi = 1$. Then an object in $\vect_\Cb^{\Zb_2}$ is a $\Cb$-vector space $V \in \vect_\Cb$ equipped with a $\Cb$-linear map $J \colon \overline{V} \to V$ such that $J \circ \bar J = 1_V$. In other words, $J \colon V \to V$ is an anti-linear map and $J^2 = 1_V$. Define
\[
V_\pm \coloneqq \{v \in V \mid J(v) = \pm v\} .
\]
Then $V_\pm$ are $\Rb$-linear subspaces of $V$ such that $V = V_+ \oplus V_-$, and $\mathrm i \cdot V_+ = V_-$. So $V = V_+ \oplus \mathrm i \cdot V_+ \simeq \Cb \otimes_\Rb V_+$. We say that $V_+ = V^{\Zb_2}$ is the \emph{real form} of $V$ determined by $J$. A morphism $f \colon (V,J) \to (W,J')$ in $\vect_\Cb^{\Zb_2}$ is a $\Cb$-linear map $f \colon V \to W$ such that $J' \circ f = f \circ J$. Therefore, it restricts to an $\Rb$-linear map $f_+ \colon V_+ \to W_+$. So taking real forms defines an $\Rb$-linear functor $\vect_\Cb^{\Zb_2} \to \vect_\Rb$, and clearly it is essentially surjective. Since $V_- = \mathrm i \cdot V_+$, the $\Cb$-linear map $f$ is completely determined by the restriction $f_+$. This shows that the functor $\vect_\Cb^{\Zb_2} \to \vect_\Rb$ is fully faithful. Moreover, it is not hard to see that the equivalence $\vect_\Cb^{\Zb_2} \simeq \vect_\Rb$ is symmetric monoidal.

For the action with $\phi = -1$, an object in $\vect_\Cb^{\Zb_2}$ is a $\Cb$-vector space $V \in \vect_\Cb$ equipped with an anti-linear map $J \colon V \to V$ such that $J^2 = -1_V$. Equivalently, $V$ is an $\Rb$-vector space equipped with two $\Rb$-linear maps $I,J \colon V \to V$ such that
\[
I^2 = J^2 = -1_V , \quad IJ = -JI .
\]
So $V$ is exactly a module over the quaternion algebra $\Hb \simeq \mathrm{Cl}_{0,2}(\Rb)$. It follows that $\vect_\Cb^{\Zb_2}$ is equivalent to the category $\vect_\Hb$ of finite-dimensional $\Hb$-modules.
\end{expl}

The first part of Example \ref{expl_semi-linear_equivariantization_Vec} is a reformulation of the classical Galois descent theory for vector spaces (over $\Cb/\Rb$). Given a $\Cb$-vector space $W$, a \emph{real form} is an $\Rb$-vector space $V$ equipped with a $\Cb$-linear isomorphism $\Cb \otimes_\Rb V \simeq W$. Example \ref{expl_semi-linear_equivariantization_Vec} shows that \emph{descent data} on $W$ that determine real forms are semi-linear $\Zb_2$-actions on $W$, i.e., anti-linear maps $J \colon W \to W$ satisfying $J^2 = 1_W$. These results can be summarized into the following theorem, which is essentially proved in Example \ref{expl_semi-linear_equivariantization_Vec}.

\begin{thm}
There is an $\Rb$-linear symmetric monoidal equivalence $\vect_\Cb^{\Zb_2} \simeq \vect_\Rb$ such that the following diagram commutes:
\[
\xymatrix{
\vect_\Cb^{\Zb_2} \ar@{<->}[rr]^-{\simeq} \ar[dr]_{\text{forget}} & & \vect_\Rb \ar[dl]^{\Cb \otimes_\Rb -} \\
 & \vect_\Cb
}
\]
where the $\Zb_2$-action on $\vect_\Cb$ is defined by the complex conjugation of $\Cb$-vector spaces with $\phi = 1$.
\end{thm}

The Galois descent theory for (tensor) categories was developed in \cite{EG11}. In this paper we only need a special case, which is formulated as the following theorem.

\begin{thm} \label{thm_Galois_descent_category}
There is an $\Rb$-linear symmetric monoidal equivalence $2\vect_\Cb^{\Zb_2} \simeq 2\vect_\Rb$ such that the following diagram commutes:
\[
\xymatrix{
2\vect_\Cb^{\Zb_2} \ar@{<->}[rr]^-{\simeq} \ar[dr]_{\text{forget}} & & 2\vect_\Rb \ar[dl]^{\vect_\Cb \boxtimes_\Rb -} \\
 & 2\vect_\Cb
}
\]
where the $\Zb_2$-action on $2\vect_\Cb$ is defined by complex conjugation of (finite semisimple) $\Cb$-linear categories. More precisely, the functor $2\vect_\Cb^{\Zb_2} \to 2\vect_\Rb$ is given by equivariantization, and the functor $2\vect_\Rb \to 2\vect_\Cb^{\Zb_2}$ is given by complexification.
\end{thm}

In the next subsection, we give a proof of Theorem \ref{thm_Galois_descent_category}. Now we discuss a corollary of this theorem.

Given a finite $\Cb$-linear category $\CD$, a \emph{real form} is a finite $\Rb$-linear category $\CC$ equipped with a $\Cb$-linear equivalence $\CC_\Cb \simeq \CD$. For example, both $\vect_\Rb$ and $\vect_\Hb$ are real forms of $\vect_\Cb$. If in addition $\CD$ is semisimple, its real forms are also semisimple because the field extension $\Cb/\Rb$ is separable. By definition, a real form of a finite semisimple $\Cb$-linear category $\CD \in 2\vect_\Cb$ is an object in the fiber over $\CD$ along the complexification functor $2\vect_\Rb \to 2\vect_\Cb$. Then by Theorem \ref{thm_Galois_descent_category}, it is equivalent to an object in the fiber over $\CD$ along the forgetful functor $2\vect_\Cb^{\Zb_2} \to 2\vect_\Rb$, which is simply a semi-linear $\Zb_2$-action on $\CD$. In other words, descent data on $\CD$ are semi-linear $\Zb_2$-actions.

Furthermore, since the equivalence $2\vect_\Cb^{\Zb_2} \simeq 2\vect_\Rb$ in Theorem \ref{thm_Galois_descent_category} is symmetric monoidal, it preserves $E_k$-algebras. Therefore, descent data on $E_k$-multi-fusion categories over $\CC$ are $E_k$-monoidal semi-linear $\Zb_2$-actions.

\begin{rem}
Not all $E_k$-multi-fusion categories over $\Cb$ have real forms. The simplest example is the fusion category $\vect_{\Cb}[\Zb_3,\omega]$, where $\omega \in H^3(\Zb_3;\Cb^\times) \simeq \Zb_3$ is a generator. Indeed, $\overline{\vect_{\Cb}[\Zb_3,\omega]} \simeq \vect_{\Cb}[\Zb_3,\omega^{-1}]$ is not monoidally equivalent to $\vect_{\Cb}[\Zb_3,\omega]$, because every automorphism of $\Zb_3$ induces the identity map on $H^3(\Zb_3;\Cb^\times)$. So there is no anti-linear autoequivalence on $\vect_{\Cb}[\Zb_3,\omega]$.

On the other hand, symmetric fusion categories over $\Cb$ always have real forms. By Deligne's theorem (see Theorem \ref{thm_Deligne_classificaion_symmetric_fusion_category}), a symmetric fusion category over $\Cb$ is equivalent to $\srep_\Cb(G,z)$ for some finite super group $(G,z)$. An obvious real form is $\srep_\Rb(G,z)$, which corresponds to the semi-linear $\Zb_2$-action on $\srep_\Cb(G,z)$ given by complex conjugation.
\end{rem}

\subsection{Proof of the Galois descent theorem}

Let us consider the category $\BMod_{\Cb|\Cb}(\vect_R)$ of $\Cb$-$\Cb$-bimodules in $\vect_\Rb$. It is a fusion category over $\Rb$, and the monoidal structure is given by the relative tensor product $\otimes_\Cb$.

\begin{lem} \label{lem_BMod_C_C_Vec_C_Z2}
There is an $\Rb$-linear equivalence $\BMod_{\Cb|\Cb}(\vect_R) \simeq \vect_\Cb[\Zb_2]$.
\end{lem}

\pf
Let $M \in \BMod_{\Cb|\Cb}(\vect_R)$. Define two $\Rb$-linear maps $L,R \colon M \to M$ by
\[
L(m) \coloneqq \mathrm i \triangleright m , \quad R(m) \coloneqq m \triangleleft \mathrm i , \quad \forall m \in M .
\]
Clearly $LR = RL$. Then
\[
(LR)^2 = L^2 R^2 = 1_M .
\]
Thus $M$ can be decomposed as the direct sum of two $\Rb$-linear subspaces $M_0 \oplus M_1$, where
\[
M_j \coloneqq \{m \in M \mid L(m) = (-1)^j \cdot R(m)\} .
\]
Since $L$ commutes with $LR$, both $M_0$ and $M_1$ are left $\Cb$-submodules of $M$. Therefore, $M \mapsto M_0 \oplus M_1$ defines a functor $\BMod_{\Cb|\Cb}(\vect_R) \to \vect_\Cb[\Zb_2]$. Its quasi-inverse is given below. For a $\Zb_2$-graded vector space $V = V_0 \oplus V_1 \in \vect_\Cb[\Zb_2]$, define a $\Cb$-$\Cb$-bimodule structure on $V$ as follows:
\[
\lambda \triangleright v \coloneqq \lambda \cdot v , \quad v \triangleleft \lambda \coloneqq \begin{cases}
\lambda \cdot v , & v \in V_0 , \\
\bar{\lambda} \cdot v , & v \in V_1 ,
\end{cases} \quad \forall \lambda \in \Cb , \, v \in V .
\]
This defines a quasi-inverse $\vect_\Cb[\Zb_2] \to \BMod_{\Cb|\Cb}(\vect_R)$. 
\epf

The $\Cb$-linear structure on $\vect_\Cb[\Zb_2]$ can be transferred to $\BMod_{\Cb|\Cb}(\vect_R)$ via the above equivalence. However, if we choose the right $\Cb$-module structure on $M = M_0 \oplus M_1$ in the proof of Lemma \ref{lem_BMod_C_C_Vec_C_Z2}, we can obtain a different equivalence $\BMod_{\Cb|\Cb}(\vect_R) \simeq \vect_\Cb[\Zb_2]$. The $\Cb$-linear structures on $\BMod_{\Cb|\Cb}(\vect_R)$ transferred from $\vect_\Cb[\Zb_2]$ via these two equivalences are not the same. So there is no canonical $\Cb$-linear structure on $\BMod_{\Cb|\Cb}(\vect_R)$. Moreover, both $\Cb$-linear structures are not compatible with the monoidal structure.

We use $\Cb_0$ to denote the regular $\Cb$-$\Cb$-bimodule, and $\Cb_1$ to denote $\Cb$ equipped with the following $\Cb$-$\Cb$-bimodule structure:
\[
\lambda \triangleright \mu \coloneqq \lambda \mu , \quad \mu \triangleleft \lambda \coloneqq \bar \lambda \mu , \quad \forall \lambda,\mu \in \Cb .
\]
In other words, $\Cb_0$ and $\Cb_1$ are the image of the two simple objects in $\vect_\Cb[\Zb_2]$ under the equivalence defined in the proof of Lemma \ref{lem_BMod_C_C_Vec_C_Z2}. The fusion rule of $\BMod_{\Cb|\Cb}(\vect_\Rb)$ is also the same as $\vect_\Cb[\Zb_2]$:
\begin{align} \label{eq_fusion_rule_BMod_C_C}
\Cb_i \otimes_\Cb \Cb_j & \simeq \Cb_{i+j} \\
\lambda \otimes_\Cb \mu & \mapsto \begin{cases}
\lambda \mu , & i = 0 , \\
\lambda \bar \mu , & i = 1 .
\end{cases} \nonumber
\end{align}
where $i+j$ is understood modulo $2$.

\begin{lem}
The 2-category $\LMod_{\BMod_{\Cb|\Cb}(\vect_\Rb)}(2\vect_\Rb)$ is equivalent to $2\vect_\Cb^{\Zb_2}$, where the $\Zb_2$-action on $2\vect_\Cb$ is defined by complex conjugation.
\end{lem}

\pf
Let $\CM$ be a finite semisimple left $\BMod_{\Cb|\Cb}(\vect_\Rb)$-module (over $\Rb$). Since $\Omega \BMod_{\Cb|\Cb}(\vect_\Rb) \simeq \Cb$, the functor $\Cb_0 \odot - \colon \CM \to \CM$ is equivalent to a $\Cb$-linear structure on $\CM$. So we view $\CM$ as a $\Cb$-linear category. Define $J \coloneqq \Cb_1 \odot - \colon \CM \to \CM$. Then $J$ is anti-linear. The module associator is determined by the natural isomorphism
\[
\phi \coloneqq \bigl( J \circ J \simeq \Cb_0 \odot - \simeq 1_\CM \bigr) .
\]
A simple calculation shows that the following diagram commutes:
\[
\xymatrix{
(\Cb_1 \otimes_\Cb \Cb_1) \otimes_\Cb \Cb_1 \ar[rr]^-{\alpha_{1,1,1}} \ar[d]_{\simeq} & & \Cb_1 \otimes_\Cb (\Cb_1 \otimes_\Cb \Cb_1) \ar[d]^{\simeq} \\
\Cb_0 \otimes_\Cb \Cb_1 \ar[r]^-{\simeq} & \Cb_1 & \Cb_1 \otimes_\Cb \Cb_0 \ar[l]_-{\simeq}
}
\]
where $\alpha$ is the associator of $\BMod_{\Cb|\Cb}(\vect_\Rb)$. It follows that the pentagon equation for the module associator is equivalent to $J(\phi_x) = \phi_{J(x)}$ for all $x \in \CM$. Hence, a left $\BMod_{\Cb|\Cb}(\vect_\Rb)$-module structure on $\CM$ is equivalent to a semi-linear $\Zb_2$-action. Similarly, one can easily check that left $\BMod_{\Cb|\Cb}(\vect_\Rb)$-module functors are $\Zb_2$-equivariant functors, and left $\BMod_{\Cb|\Cb}(\vect_\Rb)$-module natural transformations are $\Zb_2$-equivariant natural transformations.
\epf

Clearly $\vect_\Cb$ is a left $\BMod_{\Cb|\Cb}(\vect_\Rb)$-module, and the module structure is given by the relative tensor product $\otimes_\Cb$. The corresponding semi-linear $\Zb_2$-action is given by complex conjugation with $\phi = 1$. By a computation, one can show that the internal hom algebra $[\Cb,\Cb] \in \BMod_{\Cb|\Cb}(\vect_\Rb)$ is isomorphic to the matrix algebra $\rmM_2(\Rb)$, on which the $\Cb$-$\Cb$-bimodule structure is induced by the ring homomorphism
\begin{align*}
\Cb & \to \rmM_2(\Rb) \\
a+b \mathrm i & \mapsto \begin{pmatrix}
a & -b \\
b & a
\end{pmatrix}
\end{align*}
Then by Theorem \ref{thm_Ostrik_module_algebra}, $\vect_\Cb$ is equivalent to $\RMod_{\rmM_2(\Rb)}(\BMod_{\Cb|\Cb}(\vect_\Rb))$.

\begin{prop} \label{prop_semi-linear_equivariantization_relative_tensor}
Let $\CD \in 2\vect_\Cb^{\Zb_2}$. The following $\Rb$-linear categories are equivalent:
\bnu[(1)]
\item the equivariantization $\CD^{\Zb_2}$;
\item the category $\LMod_{\Rb[\Zb_2]}(\CD)$ of left $\Rb[\Zb_2]$-modules in $\CD$, where the group algebra $\Rb[\Zb_2]$ is viewed as an algebra in $\vect_\Rb[\Zb_2]$;
\item the category $\LMod_{\rmM_2(\Rb)}(\CD)$ of left $\rmM_2(\Rb)$-modules in $\CD$, where the matrix algebra $\rmM_2(\Rb)$ is viewed as an algebra in $\BMod_{\Cb|\Cb}(\vect_\Rb)$;
\item the relative Deligne tensor product $\vect_\Cb \boxtimes_{\BMod_{\Cb|\Cb}(\vect_\Rb)} \CD$.
\enu
In particular, $\CD^{\Zb_2}$ is finite semisimple.
\end{prop}

\pf
The equivalence between (1) and (2) is proved in Lemma \ref{lem_equivariantization_module_group_algebra}. Note that the canonical $\Rb$-linear monoidal functor $F \colon \vect_\Rb[\Zb_2] \to \BMod_{\Cb|\Cb}(\vect_\Rb)$ sends the group algebra $\Rb[\Zb_2] \in \vect_\Rb[\Zb_2]$ to $\rmM_2(\Rb) \in \BMod_{\Cb|\Cb}(\vect_\Rb)$. Then
\[
\LMod_{\rmM_2(\Rb)}(\BMod_{\Cb|\Cb}(\vect_\Rb)) \simeq \LMod_{F(\Rb[\Zb_2])}(\BMod_{\Cb|\Cb}(\vect_\Rb)) \simeq \LMod_{\Rb[\Zb_2]}(\BMod_{\Cb|\Cb}(\vect_\Rb)) .
\]
This proves the equivalence between (2) and (3). The equivalence between (3) and (4) is due to the equivalence $\RMod_{\rmM_2(\Rb)}(\BMod_{\Cb|\Cb}(\vect_\Rb)) \simeq \vect_\Cb$. Finally, by Theorem \ref{thm_separable_module_relative_tensor} we see that $\CD^{\Zb_2}$ is also a finite semisimple category.
\epf

Now we give a proof of Theorem \ref{thm_Galois_descent_category}.

\pf[Proof of Theorem \ref{thm_Galois_descent_category}]
The complexification functor $\vect_\Cb \boxtimes_\Rb - \colon 2\vect_\Rb \to 2\vect_\Cb$ is symmetric monoidal and can be lifted to an $\Rb$-linear 2-functor $2\vect_\Rb \to 2\vect_\Cb^{\Zb_2}$. So this lifting is also symmetric monoidal. It suffices to show that this 2-functor is an equivalence.

By Theorem \ref{thm_invertible_bimodule}, $\BMod_{\Cb|\Cb}(\vect_\Rb)$ is Morita equivalent to $\vect_\Rb$ with the invertible bimodule $\vect_\Cb$. Then by Theorem \ref{thm_Morita_equivalent_defined_by_2-category}, the $\Rb$-linear 2-functor
\[
2\vect_\Rb \simeq \LMod_{\vect_\Rb}(2\vect_\Rb) \xrightarrow{\vect_\Cb \boxtimes_{\vect_\Rb} -} \LMod_{\BMod_{\Cb|\Cb}(\vect_\Rb)}(2\vect_\Rb) \simeq 2\vect_\Cb^{\Zb_2}
\]
is an equivalence, which is the above lifting 2-functor. Its quasi-inverse is given by
\[
\vect_\Cb \boxtimes_{\BMod_{\Cb|\Cb}(\vect_\Rb)}(2\vect_\Rb) - \colon \LMod_{\BMod_{\Cb|\Cb}(\vect_\Rb)}(2\vect_\Rb) \to 2\vect_\Rb .
\]
By Proposition \ref{prop_semi-linear_equivariantization_relative_tensor}, this is the equivariantization 2-functor $(\cdot)^{\Zb_2} \colon 2\vect_\Cb^{\Zb_2} \to 2\vect_\Rb$.
\epf



\section{Representations of \texorpdfstring{$\Zb_2$}{Z2}-graded groups} \label{sec_semi-linear_rep}

\subsection{Semi-linear representations of \texorpdfstring{$\Zb_2$}{Z2}-graded groups}

\begin{defn}
A \emph{$\Zb_2$-graded group} is a group $G$ equipped with a group homomorphism $s \colon G \to \Zb_2$. A \emph{homomorphism} $f \colon (G,s) \to (G',s')$ between two $\Zb_2$-graded groups is a group homomorphism $f \colon G \to G'$ such that $s' \circ f = s$.
\end{defn}


\begin{expl}
Let $V$ be a $\Cb$-vector space. A \emph{semi-linear automorphism} of $V$ is an $\Rb$-linear map $V \to V$ that is either $\Cb$-linear or anti-linear. All semi-linear automorphisms of $V$ form a group, denoted by $\mathrm{GL}_{\Cb/\Rb}(V)$. It is equipped with an obvious group homomorphism $\mathrm{GL}_{\Cb/\Rb}(V) \to \Zb_2$ sending $\Cb$-linear maps to $0$ and anti-linear maps to $1$. So $\mathrm{GL}_{\Cb/\Rb}(V)$ is a $\Zb_2$-graded group.
\end{expl}

Let $(G,s)$ be a $\Zb_2$-graded group.

\begin{defn}
A \emph{semi-linear representation} of $(G,s)$ is a $\Cb$-vector space $V$ equipped with a $\Zb_2$-graded group homomorphism $\rho \colon (G,s) \to \mathrm{GL}_{\Cb/\Rb}(V)$. In other words, $(V,\rho)$ is an $\Rb$-linear representation such that
\[
\rho(g) \colon V \to V \text{ is} \begin{cases}
\text{$\Cb$-linear} , &  s(g) = 0 , \\
\text{anti-linear} , & s(g) = 1 ,
\end{cases} \quad \forall g \in G .
\]
A \emph{homomorphism} between two semi-linear representations $(V,\rho)$ and $(W,\sigma)$ is a $\Cb$-linear map $f \colon V \to W$ such that
\begin{equation} \label{eq_rep_morphism_condition}
f \circ \rho(g) = \sigma(g) \circ f , \quad \forall g \in G .
\end{equation}
The category of finite-dimensional semi-linear representations of $(G,s)$ is denoted by $\rep_{\Cb/\Rb}(G,s)$.
\end{defn}

It is not hard to check that $\rep_{\Cb/\Rb}(G,s)$ is an abelian category.

\begin{rem}
Although the morphisms in $\rep_{\Cb/\Rb}(G,s)$ are $\Cb$-linear maps, $\rep_{\Cb/\Rb}(G,s)$ is only an $\Rb$-linear category in general. Indeed, for $\lambda \in \Cb$, the condition \eqref{eq_rep_morphism_condition} implies that
\[
(\lambda f) \circ \rho(g) = \lambda (\sigma(g) \circ f) =  \sigma(g) \circ (\bar{\lambda} f) , \quad \text{if } s(g) = 1 .
\]
Therefore, when $f \neq 0$ and $s \colon G \to \Zb_2$ is nontrivial, $\lambda f$ is a morphism in $\rep_{\Cb/\Rb}(G,s)$ only if $\lambda \in \Rb$.
\end{rem}

Given two semi-linear representations $(V,\rho)$ and $(W,\sigma)$ of $(G,s)$, the tensor product $V \otimes_\Cb W$ is naturally a semi-linear representation defined by
\[
G \ni g \mapsto \rho(g) \otimes_\Cb \sigma(g) \colon V \otimes_\Cb W \to V \otimes_\Cb W .
\]
This defines a monoidal structure on $\rep_{\Cb/\Rb}(G,s)$. Moreover, $\rep_{\Cb/\Rb}(G,s)$ admits a symmetric monoidal structure inherited from $\vect_\Cb$, and the forgetful functor $\forget \colon \rep_{\Cb/\Rb}(G,s) \to \vect_\Cb$ is clearly an $\Rb$-linear symmetric monoidal functor.

\begin{expl}
The \emph{trivial semi-linear representation} is $\Cb$ equipped with the obvious semi-linear $(G,s)$-action $G \xrightarrow{s} \Zb_2 \simeq \Gal(\Cb/\Rb) \to \mathrm{GL}_{\Cb/\Rb}(\Cb)$. It is the tensor unit of $\rep_{\Cb/\Rb}(G,s)$.
\end{expl}

\begin{expl} \label{expl_dual_semi-linear_rep}
Let $(V,\rho) \in \rep_{\Cb/\Rb}(G,s)$. Its \emph{dual semi-linear representation} is $V^* \coloneqq \Hom_\Cb(V,\Cb)$ equipped with the semi-linear $(G,s)$-action $\rho^* \colon (G,s) \to \mathrm{GL}_{\Cb/\Rb}(V^*)$ defined by
\[
\rho^*(g)(\alpha) \coloneqq \overline{(\cdot)}^{s(g)} \circ \alpha \circ \rho(g)^{-1} , \quad g \in G , \, \alpha \in V^* .
\]
It is both the left and right dual of $(V,\rho)$ with the usual evaluation and coevaluation maps.
\end{expl}

\begin{expl} \label{expl_regular_rep}
Suppose $G$ is finite. Let $\fun^l(G,s)$ be the function space $\fun(G,\Cb)$ equipped with the left translation semi-linear $(G,s)$-action defined as follows:
\[
(g \triangleright f)(x) \coloneqq \begin{cases}
f(g^{-1}x) , & s(g) = 0 , \\
\overline{f(g^{-1}x)} , & s(g) = 1 ,
\end{cases} \quad g,x \in G , \, f \in \fun(G,\Cb) .
\]
This is called the \emph{left regular semi-linear representation} of $(G,s)$. Similarly, $\fun^r(G,s)$ is the function space $\fun(G,\Cb)$ equipped with the right translation semi-linear $(G,s)$-action defined as follows:
\[
(g \triangleright f)(x) \coloneqq \begin{cases}
f(xg) , & s(g) = 0 , \\
\overline{f(xg)} , & s(g) = 1 ,
\end{cases} \quad g,x \in G , \, f \in \fun(G,\Cb) .
\]
This is called the \emph{right regular semi-linear representation} of $(G,s)$. They are isomorphic via $f \mapsto f^\vee$, where $f^\vee(x) \coloneqq f(x^{-1})$ for all $x \in G$. Moreover, they are commutative algebras in $\rep_{\Cb/\Rb}(G,s)$ with the pointwise multiplication and are called the \emph{regular algebras}.
\end{expl}

\begin{rem}
Let $\forget \colon \rep_{\Cb/\Rb}(G,s) \to \vect_\Cb$ be the forgetful functor and $\forget^R$ be the right adjoint. Then $\forget^R(\Cb) \in \rep_{\Cb/\Rb}(G,s)$ is a commutative separable algebra \cite[Lemma 3.5]{DMNO13}. It is not hard to check that $\forget^R(\Cb)$ is isomorphic to the regular algebra defined in Example \ref{expl_regular_rep}.
\end{rem}

\begin{expl}
When $s \colon G \to \Zb_2$ is trivial, $\rep_{\Cb/\Rb}(G,s)$ is the same as $\rep_\Cb(G)$ as symmetric monoidal categories.
\end{expl}

\begin{expl} \label{expl_Z2T_rep}
We use $\Zb_2^T$ to denote the $\Zb_2$-graded group $(\Zb_2,1 \colon \Zb_2 \to \Zb_2)$. Here the superscript ``$T$'' stands for ``time-reversal''. A semi-linear representation of $\Zb_2^T = (\Zb_2,1)$ is a $\Cb$-linear vector space $V$ equipped with an anti-linear map $J \colon V \to V$ such that $J^2 = 1_V$. By Example \ref{expl_semi-linear_equivariantization_Vec}, $\rep_{\Cb/\Rb}(\Zb_2^T) \simeq \vect_\Cb^{\Zb_2} \simeq \vect_\Rb$ as symmetric monoidal categories.
\end{expl}

\begin{expl} \label{expl_Z2T_projective_rep}
We are also interested in the projective semi-linear representations of $\Zb_2^T$. A \emph{projective semi-linear representation} of $\Zb_2^T$ is a $\Cb$-linear vector space $V$ equipped with an anti-linear map $J \colon V \to V$ such that $J^2 = \lambda \cdot 1_V$ for some $\lambda \in \Cb^\times$. Note that
\[
\lambda \cdot J = (JJ)J = J(JJ) = \bar{\lambda} \cdot J .
\]
Thus $\lambda \in \Rb^\times$. We can further scale $J$ by a real number so that $\lambda = \pm 1$. In quantum mechanics, these two kinds of projective semi-linear representations are called the time-reversal symmetries with $T^2 = \pm 1$, respectively.

When $\lambda = 1$, this is discussed in Example \ref{expl_Z2T_rep}. When $\lambda = -1$, by Example \ref{expl_semi-linear_equivariantization_Vec}, a projective semi-linear representation is exactly an $\Hb$-module. So the category of finite-dimensional projective semi-linear representations of $\Zb_2^T$ with $T^2 = -1$ is equivalent to the category $\vect_\Hb$ of finite-dimensional $\Hb$-modules.
\end{expl}

\begin{expl}
Let $G = K \times \Zb_2$ where $K$ is a group, and $s \colon K \times \Zb_2 \to \Zb_2$ be the projection to the second factor. We also denote $(K \times \Zb_2,s)$ by $K \times \Zb_2^T$. Similar to Example \ref{expl_Z2T_rep}, $\rep_{\Cb/\Rb}(K \times \Zb_2^T)$ is equivalent to $\rep_\Rb(K)$ as symmetric monoidal categories.
\end{expl}

\begin{expl} \label{expl_semi-linear_rep_semidirect_product}
Let $G = K \rtimes \Zb_2$ where $K$ is a group equipped with a $\Zb_2$-action, and $s \colon K \rtimes \Zb_2 \to \Zb_2$ be the projection to the second factor. We also denote $(K \rtimes \Zb_2,s)$ by $K \rtimes \Zb_2^T$. The embedding $\Zb_2 \to K \rtimes \Zb_2$ preserves the $\Zb_2$-grading, and thus induces an $\Rb$-linear fiber functor $\rep_{\Cb/\Rb}(K \rtimes \Zb_2^T) \to \rep_{\Cb/\Rb}(\Zb_2^T) \simeq \vect_\Rb$.

Suppose $K$ is a finite abelian group and the $\Zb_2$-action on $K$ is given by $k \mapsto k^{-1}$. Then every $(V,\rho) \in \rep_{\Cb/\Rb}(K \rtimes \Zb_2^T)$, as a $\Cb$-linear representation of $K$, is the direct sum of invariant subspaces
\[
V = \bigoplus_{\chi \in \hat K} V_\chi ,
\]
where $\hat K = \Hom(K,\Cb^\times)$ is the dual group of $K$ and
\[
V_\chi \coloneqq \{v \in V \mid \rho(k,0)(v) = \chi(k) \cdot v , \, \forall k \in K \} .
\]
Denote $J \coloneqq \rho(e,1)$. Then for every $v \in V_\chi$, we have
\[
\rho(k,0)(J(v)) = J(\rho(k^{-1},0)(v)) = J(\chi(k^{-1}) \cdot v) = \chi(k) \cdot J(v) , \quad \forall k \in K .
\]
In other words, $V_\chi$ is invariant under $J$. Then by Example \ref{expl_Z2T_rep}, $\rep_{\Cb/\Rb}(K \rtimes \Zb_2^T)$ is equivalent to $\vect_\Rb[\hat K]$ as $\Rb$-linear monoidal categories. Clearly the braiding structure is trivial.
\end{expl}

\begin{expl} \label{expl_Z4T_rep}
Let $G = \Zb_4$ and $s \colon \Zb_4 \to \Zb_2$ be the nontrivial group homomorphism (i.e., modulo $2$). We also denote $(\Zb_4,s)$ by $\Zb_4^T$. A semi-linear representation of $\Zb_4^T$ is a $\Cb$-linear vector space $V$ equipped with an anti-linear map $J \colon V \to V$ such that $J^4 = 1_V$. Then $J^2$ is a $\Cb$-linear map on $V$ such that $(J^2)^2 = 1_V$. Define
\[
V_\pm \coloneqq \{v \in V \mid J^2 (v) = \pm v\} .
\]
Then $V_\pm$ are $\Cb$-linear subspaces of $V$ and $V = V_+ \oplus V_-$. For $v \in V_\pm$, we have
\[
J^2(J(v)) = J(J^2(v)) = J(\pm v) = \pm J(v) .
\]
Thus $V_\pm$ are also invariant under $J$. It follows that $(V_\pm,J\vert_{V_\pm})$ are projective semi-linear representations of $\Zb_2^T$ with $T^2 = \pm 1$ (see Example \ref{expl_Z2T_projective_rep}). By Examples \ref{expl_Z2T_rep} and \ref{expl_Z2T_projective_rep}, $\rep_{\Cb/\Rb}(\Zb_4^T)$ is equivalent to $\vect_\Rb \oplus \vect_\Hb$ as $\Rb$-linear categories.

We also explicitly write down the simple objects and fusion rule of $\rep_{\Cb/\Rb}(\Zb_4^T)$. One is $\Cb$ on which the generator of $\Zb_4^T$ acts as complex conjugation. The other one is $\Hb$ on which the generator of $\Zb_4^T$ acts as $J$. The fusion rule is
\[
\Hb \otimes_\Cb \Hb \simeq \mathrm M_2(\Cb) \simeq \Cb^{\oplus 4}
\]
and $\Cb$ is the tensor unit.
\end{expl}


\begin{defn}
A \emph{semi-unitary representation} of $(G,s)$ is a semi-linear representation $(V,\rho)$ equipped with a Hilbert space structure such that $\rho(g)$ is unitary if $s(g) = 0$ and is anti-unitary if $s(g) = 1$.
\end{defn}

\begin{prop}
Every finite-dimensional semi-linear representation of $(G,s)$ admits a positive definite Hermitian inner product that makes it into a semi-unitary representation.
\end{prop}

\pf
The proof is standard. Choose any positive definite Hermitian inner product $\langle -,- \rangle$ on $V$. Define
\[
\langle v,w \rangle_G \coloneqq \sum_{\substack{g \in G\\s(g) = 0}} \langle \rho(g)(v),\rho(g)(w) \rangle + \sum_{\substack{g \in G\\s(g) = 1}} \langle \rho(g)(w),\rho(g)(v) \rangle .
\]
Then $\langle -,- \rangle_G$ is also a positive definite Hermitian inner product. Clearly $(V,\rho)$ equipped with $\langle -,- \rangle_G$ is a semi-unitary representation.
\epf

\subsection{Twisted group algebras of \texorpdfstring{$\Zb_2$}{Z2}-graded groups}

Let $(G,s)$ be a $\Zb_2$-graded finite group.

\begin{defn}[\text{\cite[Section 6.2]{KZ24}}]
The \emph{twisted group algebra} $\Cb[G,s]$ of $(G,s)$ is the $\Rb$-algebra defined as follows.
\bit
\item As a vector space, $\Cb[G,s] = \Cb[G]$. Thus it has an $\Rb$-basis
\[
\{g \mid g \in G\} \sqcup \{\mathrm i g \mid g \in G \} .
\]
\item The multiplication is given by
\[
g \cdot h \coloneqq gh , \quad (\mathrm i g) \cdot h \coloneqq \mathrm i (gh) , \quad g \cdot (\mathrm i h) \coloneqq \begin{cases}
\mathrm i (gh) , & s(g) = 0 , \\
-\mathrm i (gh) , & s(g) = 1 ,
\end{cases} \quad
(\mathrm i g) \cdot (\mathrm i h) \coloneqq \begin{cases}
- gh , & s(g) = 0 , \\
gh , & s(g) = 1 .
\end{cases}
\]
In other words, $g \cdot \mathrm i = (-1)^{s(g)} \mathrm i \cdot g$.
\item The unit is given by the unit of $G$.
\eit
\end{defn}

Clearly $\Cb[G,s]$ is a semi-linear representation of $(G,s)$, where the $G$-action is given by left multiplication. It is isomorphic to the regular semi-linear representation defined in Example \ref{expl_regular_rep}.

\begin{lem}
Let $(V,\rho) \in \rep_{\Cb/\Rb}(G,s)$. There is a natural isomorphism $\rep_{\Cb/\Rb}(G,s)(\Cb[G,s],(V,\rho)) \simeq V$. Thus $\Cb[G,s] \otimes_\Rb - \colon \vect_\Rb \to \rep_{\Cb/\Rb}(G,s)$ is left adjoint to the forgetful functor $\rep_{\Cb/\Rb}(G,s) \to \vect_\Rb$.
\end{lem}

\pf
For any $v \in V$, define a $\Cb$-linear map $f_v \colon \Cb[G,s] \to V$ by
\[
f_v(g) \coloneqq \rho(g)(v) , \quad \forall g \in G .
\]
This is a homomorphism of semi-linear representations because
\[
f_v(gh) = \rho(gh)(v) = \rho(g)(\rho(h)(v)) = \rho(g)(f_v(h)) , \quad \forall g,h \in G .
\]
Then $v \mapsto f_v$ defines a map $V \to \rep_{\Cb/\Rb}(G,s)(\Cb[G,s],(V,\rho))$. It is injective because $f_v(e) = v$ for all $v \in V$. 

Let $f \colon \Cb[G,s] \to (V,\rho)$ be a homomorphism of semi-linear representations. Denote $v \coloneqq f(e)$. Then
\[
f(g) = f(ge) = \rho(g)(f(e)) = \rho(g)(v) , \quad \forall g \in G .
\]
So $f = f_v$. This shows that $v \mapsto f_v$ is surjective.
\epf

\begin{prop} \label{prop_time-reversal_rep_twisted_group_algebra}
There is an equivalence $\rep_{\Cb/\Rb}(G,s) \simeq \LMod_{\Cb[G,s]}(\vect_\Rb)$ of $\Rb$-linear categories.
\end{prop}

\pf
We already know that $\rep_{\Cb/\Rb}(G,s)$ is an abelian category. The forgetful functor $\forget \colon \rep_{\Cb/\Rb}(G,s) \to \vect_\Rb$ is clearly faithful and exact. Then by the monadicity theorem, the adjunction $(\Cb[G,s] \otimes_\Rb -) \dashv \forget$ is monadic. Therefore, $\rep_{\Cb/\Rb}(G,s)$ is equivalent to the Eilenberg-Moore category of the induced monad $\Cb[G,s] \otimes_\Rb -$ on $\vect_\Rb$. This monad structure is equivalent to an algebra structure on the space $\Cb[G,s]$: the multiplication map $\Cb[G,s] \otimes_\Rb \Cb[G,s] \to \Cb[G,s]$ is the mate of the identity morphism on $\Cb[G,s]$ under the adjunction. It is not hard to check that this algebra is nothing but the twisted group algebra $\Cb[G,s]$. 
\epf

\begin{lem} \label{lem_twisted_group_algebra_separable}
The twisted group algebra $\Cb[G,s]$ is separable.
\end{lem}

\pf
Define an $\Rb$-linear map $\delta \colon \Cb[G,s] \to \Cb[G,s] \otimes_\Rb \Cb[G,s]$ by
\begin{align*}
\delta(g) & \coloneqq \frac{1}{2\lvert G \rvert} \biggl(\sum_{h \in G} gh \otimes_\Rb h^{-1} - \sum_{h \in G} (-1)^{s(gh)} \cdot (\mathrm i gh) \otimes_\Rb (\mathrm i h^{-1}) \biggr) , \\
\delta(\mathrm i g) & \coloneqq \frac{1}{2\lvert G \rvert} \biggl(\sum_{h \in G} (\mathrm i gh) \otimes_\Rb h^{-1} + \sum_{h \in G} (-1)^{s(gh)} \cdot gh \otimes_\Rb (\mathrm i h^{-1}) \biggr) .
\end{align*}
One can check that $\delta$ is a $\Cb[G,s]$-$\Cb[G,s]$-bimodule map and $\mu \circ \delta = 1$ where $\mu$ is the multiplication of $\Cb[G,s]$.
\epf

\begin{rem}
Moreover, $\Cb[G,s]$ is a separable Frobenius algebra with the counit $\varepsilon(e) = \varepsilon(\mathrm i e) = 1$ and $\varepsilon(g) = \varepsilon(\mathrm i g) = 0$ for all $g \neq e$.
\end{rem}

\begin{prop}
$\rep_{\Cb/\Rb}(G,s)$ is a symmetric fusion category over $\Rb$.
\end{prop}

\pf
By Proposition \ref{prop_time-reversal_rep_twisted_group_algebra}, $\rep_{\Cb/\Rb}(G,s)$ is an $\Rb$-linear finite category. By Lemma \ref{lem_twisted_group_algebra_separable}, $\rep_{\Cb/\Rb}(G,s)$ is semisimple. The tensor product $\otimes_\Cb$ on $\rep_{\Cb/\Rb}(G,s)$ is exact in each variable. The dual of an object $(V,\rho) \in \rep_{\Cb/\Rb}(G,s)$ is given by the dual semi-linear representation defined in Example \ref{expl_dual_semi-linear_rep}.
\epf

\begin{rem}
The twisted group algebra $\Cb[G,s]$ is the crossed product (semidirect product) $G \rtimes \Cb$, where $G$ acts on $\Cb$ by complex conjugation via
\[
G \xrightarrow{s} \Zb_2 \simeq \Gal(\Cb/\Rb) .
\]
Then $\Cb$ equipped with this $G$-action is an algebra in the fusion category $\rep_\Rb(G)$. There is an equivalence $\LMod_{\Cb}(\rep_\Rb(G)) \simeq \LMod_{\Cb[G,s]}(\vect_\Rb) \simeq \rep_{\Cb/\Rb}(G,s)$ of $\Rb$-linear categories. Moreover, under this equivalence, the forgetful functor $\LMod_\Cb(\rep_\Rb(G)) \to \rep_\Rb(G)$ can be identified with the forgetful functor $\rep_{\Cb/\Rb}(G,s) \to \rep_\Rb(G)$.
\end{rem}

\begin{expl}
The twisted group algebra $\Cb[\Zb_2^T]$ is the $\Rb$-algebra with two generators $\mathrm i$ and $\dagger$ subject to the relations
\[
\mathrm i^2 = -1 , \quad \dagger^2 = 1 , \quad \dagger \cdot \mathrm i = -\mathrm i \cdot \dagger .
\]
Therefore, $\Cb[\Zb_2^T]$ is isomorphic to the Clifford algebra $\mathrm{Cl}_{1,1}(\Rb) \simeq \mathrm M_2(\Rb)$, which is Morita equivalent to $\Rb$. Hence we obtain the equivalence $\rep_{\Cb/\Rb}(\Zb_2^T) \simeq \LMod_{\Cb[\Zb_2^T]}(\vect_\Rb) \simeq \vect_\Rb$ again.
\end{expl}

\begin{rem}
The twisted group algebra $\Cb[G,s]$ admits a weak Hopf algebra structure defined as follows:
\bit
\item The comultiplication is
\[
\Delta(g) \coloneqq \frac12 (g \otimes_\Rb g - \mathrm i g \otimes_\Rb \mathrm i g) , \quad \Delta(\mathrm i g) \coloneqq \frac12 (g \otimes_\Rb \mathrm i g + \mathrm i g \otimes_\Rb g) .
\]
\item The counit is
\[
\varepsilon(g) \coloneqq 2 , \quad \varepsilon(\mathrm i g) \coloneqq 0 .
\]
\item The antipode map is
\[
S(g) \coloneqq g^{-1} , \quad S(\mathrm i g) \coloneqq g^{-1} \mathrm i = (-1)^{s(g)} \cdot \mathrm i g^{-1} .
\]
\eit
Via the Tannaka-Krein reconstruction theorem for weak Hopf algebras \cite{Szl01,Szl05}, the above weak Hopf algebra structure can be reconstructed from $\rep_{\Cb/\Rb}(G,s) \to \vect_\Cb \to \vect_\Rb$, which is a separable Frobenius monoidal functor.
\end{rem}


\subsection{Semi-linear super representations of \texorpdfstring{$\Zb_2$}{Z2}-graded super groups}

\begin{defn}
A \emph{$\Zb_2$-graded super group} is a super group $(G,z)$ equipped with a group homomorphism $s \colon G \to \Zb_2$ such that $s(z) = 0$. A \emph{homomorphism} $f \colon (G,z,s) \to (G',z',s')$ between two $\Zb_2$-graded super groups is a group homomorphism $f \colon G \to G'$ such that $f(z) = z'$ and $s' \circ f = s$.
\end{defn}

\begin{expl}
Let $V = V_0 \oplus V_1$ be a $\Zb_2$-graded $\Cb$-vector space (super $\Cb$-vector space). The group $\mathrm{GL}_{\Cb/\Rb}(V)_0$ of grading-preserving semi-linear automorphisms of $V$ together with the parity automorphism and the obvious group homomorphism $\mathrm{GL}_{\Cb/\Rb}(V)_0 \to \Zb_2$ is a $\Zb_2$-graded super group.
\end{expl}

Let $(G,z,s)$ be a $\Zb_2$-graded super group.

\begin{defn}
A \emph{semi-linear super representation} of $(G,z,s)$ is a super $\Cb$-vector space $V = V_0 \oplus V_1$ equipped with a $\Zb_2$-graded super group homomorphism $\rho \colon (G,z,s) \to \mathrm{GL}_{\Cb/\Rb}(V)_0$. In other words, $(V,\rho)$ is a semi-linear representation of $(G,s)$ such that $\rho(z) = \Pi_V$ is the parity automorphism. A \emph{homomorphism} between two semi-linear super representations $(V,\rho)$ and $(W,\sigma)$ is a grading-preserving $\Cb$-linear map $f \colon V \to W$ such that
\[
f \circ \rho(g) = \sigma(g) \circ f , \quad \forall g \in G .
\]
The category of finite-dimensional semi-linear super representations of $(G,z,s)$ is denoted by $\srep_{\Cb/\Rb}(G,z,s)$.
\end{defn}

\begin{lem}
The functor $\srep_{\Cb/\Rb}(G,z,s) \to \rep_{\Cb/\Rb}(G,s)$ defined by forgetting the grading is an equivalence.
\end{lem}

\pf
By definition it is fully faithful. Now we show that it is essentially surjective. Let $(V,\rho) \in \rep_{\Cb/\Rb}(G,s)$. Then $\rho(z)$ is a $\Cb$-linear automorphism of $V$ satisfying $\rho(z)^2 = 1_V$. So $V$ can be decomposed into $V = V_0 \oplus V_1$, where
\[
V_i = \{v \in V \mid \rho(z)(v) = (-1)^i \cdot v\}
\]
are the eigenspaces of $\rho(z)$. It is easy to see that $V = V_0 \oplus V_1$ is a semi-linear super representation of $(G,z,s)$.
\epf

Therefore, when $G$ is finite, $\srep_{\Cb/\Rb}(G,z,s)$ is an $\Rb$-linear fusion category. It admits a symmetric braiding inherited from $\svect_\Cb$, and the forgetful functor $\srep_{\Cb/\Rb}(G,z,s) \to \svect_\Cb$ is an $\Rb$-linear fiber functor.

\begin{expl}
When $z = e$ is trivial, $\srep_{\Cb/\Rb}(G,z,s)$ is the same as $\rep_{\Cb/\Rb}(G,s)$ as symmetric monoidal categories. When $s \colon G \to \Zb_2$ is trivial, $\srep_{\Cb/\Rb}(G,z,s)$ is the same as $\srep_\Cb(G,z)$ as symmetric monoidal categories.
\end{expl}

We use $\Zb_2^f$ to denote the super group $(\Zb_2,1)$. Here the superscript ``f'' stands for ``fermionic''.

\begin{expl}
Let $G = K \times \Zb_2$ where $(K,z)$ is a finite super group, and $s \colon K \times \Zb_2 \to \Zb_2$ be the projection to the second factor. We denote $(K \times \Zb_2,(z,0),s)$ by $K^f \times \Zb_2^T$. Then $\srep_{\Cb/\Rb}(K^f \times \Zb_2^T)$ is equivalent to $\srep_\Rb(K^f)$ as symmetric fusion categories. In particular, $\srep_{\Cb/\Rb}(\Zb_2^f \times \Zb_2^T) \simeq \svect_\Rb$ as symmetric fusion categories.
\end{expl}

\begin{expl}
Let $(G,s) = K \rtimes \Zb_2^T$ where $(K,z)$ is a finite super group equipped with a $\Zb_2$-action such that $z \in K$ is a fixed point. Then $(K,z) \rtimes \Zb_2^T$ is a $\Zb_2$-graded finite super group. The super group homomorphism $\Zb_2^f \times \Zb_2^T \to (K,z) \rtimes \Zb_2^T$ preserves the $\Zb_2$-grading, and thus induces an $\Rb$-linear fiber functor $\srep_{\Cb/\Rb}((K,z) \rtimes \Zb_2^T) \to \srep_{\Cb/\Rb}(\Zb_2^f \times \Zb_2^T) \simeq \svect_\Rb$.
\end{expl}

\begin{expl}
Let $(G,s) = K \rtimes \Zb_2^T$ where $K$ is a finite abelian group equipped with the $\Zb_2$-action given by $k \mapsto k^{-1}$. By Example \ref{expl_semi-linear_rep_semidirect_product}, $\rep_{\Cb/\Rb}(K \rtimes \Zb_2^T) \simeq \vect_\Rb[\hat K]$ as $\Rb$-linear fusion categories. The symmetric structures on $\vect_\Rb[\hat K]$ are one-to-one corresponding to the symmetric structures on the maximal 2-group, whose first homotopy group is $\hat K$ and second homotopy group is $\Rb^\times$. Then by \cite[Theorem 3.3]{JS93}, the symmetric structures on $\vect_\Rb[\hat K]$ are classified by the group of homomorphisms $\hat K \to \Rb^\times$ of order $2$, which is isomorphic to $\prescript{}{2}{K} \coloneqq \{z \in K \mid z^2 = e\}$. For any $z \in \prescript{}{2}{K}$, the fusion category $\vect_\Rb[\hat K]$ equipped with the corresponding symmetric structure is denoted by $\vect_\Rb[\hat K,z]$. The only nontrivial braidings in $\vect_\Rb[\hat K,z]$ are given by
\begin{align*}
\Cb_\chi \otimes_\Cb \Cb_\chi & \to \Cb_\chi \otimes_\Cb \Cb_\chi \\
v \otimes_\Cb w & \mapsto \chi(z) \cdot w \otimes_\Cb v ,
\end{align*}
where $\Cb_\chi$ is the one-dimensional representation corresponding to $\chi \in \hat K$. It follows that $\srep_{\Cb/\Rb}((K,z) \rtimes \Zb_2^T)$ and $\vect_\Rb[\hat K,z]$ are equivalent as symmetric fusion categories.
\end{expl}

\begin{expl}
Let $G = \Zb_4$, $z = 2 \in \Zb_4$, and $s \colon \Zb_4 \to \Zb_2$ be the nontrivial group homomorphism (i.e., modulo $2$). We denote $(\Zb_4,z,s)$ by $\Zb_4^{f,T}$. Then $\srep_{\Cb/\Rb}(\Zb_4^{f,T})$ and $\rep_{\Cb/\Rb}(\Zb_4^T)$ are equivalent as fusion categories. By Example \ref{expl_Z4T_rep}, there are two simple objects $\Cb$ and $\Hb$. Clearly the action of $z$ on $\Cb$ is $1$ and on $\Hb$ is $-1$. So the only difference between $\srep_{\Cb/\Rb}(\Zb_4^{f,T})$ and $\rep_{\Cb/\Rb}(\Zb_4^T)$ is a minus sign in the self-braiding of $\Hb$.
\end{expl}

\subsection{Complexification}

Let $(G,s)$ be a $\Zb_2$-graded finite group. There is an obvious restriction functor $\rep_{\Cb/\Rb}(G,s) \to \rep_\Cb(\ker(s))$, which is an $\Rb$-linear symmetric monoidal functor. The following result can be found in \cite[Proposition* 6.11]{KZ24}.

\begin{prop} \label{prop_complexification_rep}
Suppose $s$ is nontrivial. The restriction functor $\rep_{\Cb/\Rb}(G,s) \to \rep_\Cb(\ker(s))$ is a complexification of $\rep_{\Cb/\Rb}(G,s)$.
\end{prop}

\pf
Denote $K \coloneqq \ker(s)$. First we show that $\Cb[G,s] \otimes_\Rb \Cb$ is Morita equivalent to $\Cb[K]$ with the invertible bimodule $\Cb[G,s]$.
\bit
\item We equip $\Cb[G,s]$ with the following $\Cb[K]$-$(\Cb[G,s] \otimes_\Rb \Cb)$-bimodule structure. The left $\Cb[K]$-module structure is induced by the inclusion $\Cb[K] \hookrightarrow \Cb[G,s]$ and the regular left $\Cb[G,s]$-module structure. The right $\Cb[G,s] \otimes_\Rb \Cb$-module structure is induced by the regular right $\Cb[G,s]$-module structure and the $\Cb$-linear structure on $\Cb[G,s]$. 
\item The space $\End_{(\Cb[G,s] \otimes_\Rb \Cb)^\op}(\Cb[G,s])$ is the subspace of $\End_{\Cb[G,s]^\op}(\Cb[G,s])$ consisting of $\Cb$-linear maps. The algebra isomorphism $\End_{\Cb[G,s]^\op}(\Cb[G,s]) \simeq \Cb[G,s]$ means that every right $\Cb[G,s]$-module endomorphism of $\Cb[G,s]$ is the left multiplication of an element $x \in \Cb[G,s]$. By definition, this map is $\Cb$-linear if and only if $x \in \Cb[K]$. So we have $\End_{(\Cb[G,s] \otimes_\Rb \Cb)^\op}(\Cb[G,s]) \simeq \Cb[K]$ as $\Cb$-linear algebras.
\eit
Therefore, by Proposition \ref{prop_complexification_algebra}, the following functor is a complexification of $\rep_{\Cb/\Rb}(G,s)$:
\begin{multline*}
\rep_{\Cb/\Rb}(G,s) \simeq \LMod_{\Cb[G,s]}(\vect_\Rb) \xrightarrow{(\Cb[G,s] \otimes_\Rb \Cb) \otimes_{\Cb[G,s]} -} \LMod_{\Cb[G,s] \otimes_\Rb \Cb}(\vect_\Rb) \\
\simeq \LMod_{\Cb[G,s] \otimes_\Rb \Cb}(\vect_\Cb) \xrightarrow{\Cb[G,s] \otimes_{(\Cb[G,s] \otimes_\Rb \Cb)} -} \LMod_{\Cb[K]}(\vect_\Cb) \simeq \rep_\Cb(K) .
\end{multline*}
This functor sends $V \in \rep_{\Cb/\Rb}(G,s)$ to
\[
\Cb[G,s] \otimes_{(\Cb[G,s] \otimes_\Rb \Cb)} (\Cb[G,s] \otimes_\Rb \Cb) \otimes_{\Cb[G,s]} V \simeq V \in \rep_\Cb(K) .
\]
Hence it is isomorphic to the restriction functor.
\epf

Similarly, for any $\Zb_2$-graded finite super group $(G,z,s)$, there is an obvious restriction functor $\srep_{\Cb/\Rb}(G,z,s) \to \srep_\Cb(\ker(s),z)$, which is an $\Rb$-linear symmetric monoidal functor.

\begin{prop}
Suppose $s$ is nontrivial. The restriction functor $\srep_{\Cb/\Rb}(G,z,s) \to \srep_\Cb(\ker(s),z)$ is a complexification of $\srep_{\Cb/\Rb}(G,z,s)$.
\end{prop}

\pf
As a monoidal functor, the restriction functor is the same as $\rep_{\Cb/\Rb}(G,s) \to \rep_\Cb(\ker(s))$, which is a complexification by Proposition \ref{prop_complexification_rep}. The statement follows from the fact that this functor is symmetric.
\epf

\section{Symmetric fusion categories over \texorpdfstring{$\Rb$}{R}} \label{sec_symmetric_fusion_R}

\subsection{Semi-linear autoequivalences}

Let $G$ be a finite group. The 2-group of semi-linear symmetric autoequivalences of $\rep_\Cb(G)$ and monoidal natural isomorphisms is denoted by $\Aut_{\Cb/\Rb}^{E_3}(\rep_\Cb(G))$. It contains $\Aut_\Cb^{E_3}(\rep_\Cb(G))$ as a monoidal full subcategory.

\begin{lem} \label{lem_auto_rep_G_semi_linear}
The 2-group $\Aut_{\Cb/\Rb}^{E_3}(\rep_\Cb(G))$ is equivalent to $\Aut_\Cb^{E_3}(\rep_\Cb(G)) \times \Zb_2$.
\end{lem}

\pf
The complex conjugation $\overline{(\cdot)} \colon \rep_\Cb(G) \to \rep_\Cb(G)$ is an anti-linear symmetric monoidal equivalence. By Proposition \ref{prop_auto_rep_G_linear}, every $F \in \Aut_\Cb^{E_3}(\rep_\Cb(G))$ is isomorphic to $\phi^*$ for some group automorphism $\phi$ of $G$. So $F \circ \overline{(\cdot)} \simeq \overline{(\cdot)} \circ F$ as symmetric monoidal functors. If $F'$ is an anti-linear symmetric monoidal autoequivalence of $\rep_\Cb(G)$, the composition $F' \circ \overline{(\cdot)}$ is $\Cb$-linear. So we have
\[
F' \circ \overline{(\cdot)} = \bigl(F' \circ \overline{(\cdot)}\bigr) \circ \overline{(\cdot)} \circ \overline{(\cdot)} \simeq \overline{(\cdot)} \circ \bigl(F' \circ \overline{(\cdot)}\bigr) \circ \overline{(\cdot)} = \overline{(\cdot)} \circ F' .
\]

Define
\begin{align*}
\Aut_{\Cb/\Rb}^{E_3}(\rep_\Cb(G)) & \to \Aut_\Cb^{E_3}(\rep_\Cb(G)) \times \Zb_2 \\
F & \mapsto \begin{cases}
(F,0) , & F \text{ is $\Cb$-linear}, \\
(F \circ \overline{(\cdot)},1) , & F \text{ is anti-linear} ,
\end{cases}
\end{align*}
and
\begin{align*}
\Aut_\Cb^{E_3}(\rep_\Cb(G)) \times \Zb_2 & \to \Aut_{\Cb/\Rb}^{E_3}(\rep_\Cb(G)) \\
(F,n) & \mapsto F \circ \overline{(\cdot)}^n .
\end{align*}
They are monoidal functors and inverse to each other.
\epf

Combining Proposition \ref{prop_auto_rep_G_linear} and Lemma \ref{lem_auto_rep_G_semi_linear}, we obtain the following result.

\begin{prop} \label{prop_auto_rep_G_semi_linear}
The 2-group $\Aut_{\Cb/\Rb}^{E_3}(\rep_\Cb(G))$ is equivalent to $\Aut(\rmB G) \times \Zb_2$.
\end{prop}


Let $(G,z)$ be a finite super group.

\begin{lem} \label{lem_auto_srep_G_z_semi_linear}
The 2-group $\Aut_{\Cb/\Rb}^{E_3}(\srep_\Cb(G,z))$ of semi-linear symmetric autoequivalences of $\srep_\Cb(G,z)$ is equivalent to $\Aut_\Cb^{E_3}(\srep_\Cb(G,z)) \times \Zb_2$.
\end{lem}

\pf
The proof is completely identical to that of Lemma \ref{lem_auto_rep_G_semi_linear}.
\epf

Then we obtain the following result.

\begin{prop} \label{prop_auto_srep_G_z_semi_linear}
The 2-group $\Aut_{\Cb/\Rb}^{E_3}(\srep_\Cb(G,z))$ is equivalent to $\Aut_{\rmB \Zb_2}(\rm B G) \times \Zb_2$.
\end{prop}

\begin{rem}
Lemmas \ref{lem_auto_rep_G_semi_linear} and  \ref{lem_auto_srep_G_z_semi_linear} are special cases of the following result. Let $\CC$ be a symmetric fusion category over $\Cb$. By \cite[Proposition 3.5]{EG11}, it has a real form if and only if the following extension of 2-groups splits:
\[
\xymatrix{
1 \ar[r] & \Aut^{E_3}_\Cb(\CC) \ar[r] & \Aut^{E_3}_{\Cb/\Rb}(\CC) \ar[r] & \Gal(\Cb/\Rb) \simeq \Zb_2 \ar[r] & 1
}
\]
In other words, $\Aut^{E_3}_{\Cb/\Rb}(\CC)$ is equivalent to $\Aut^{E_3}_\Cb(\CC) \times \Zb_2$ as 2-groups.
\end{rem}

\subsection{Group extensions and equivariantizations} \label{sec_extension_equivariantization}

Let $\Gamma,K$ be finite groups. It is well-known that monoidal functors $\Gamma \to \Aut(\rmB K)$ correspond to (non-abelian) group extensions of $K$ by $\Gamma$:
\[
\xymatrix{
1 \ar[r] & K \ar[r] & G \ar[r]^-{s} & \Gamma \ar[r] & 1
}
\]
or equivalently, $\Gamma$-graded groups $(G,s)$ such that $s$ is surjective and $\ker(s) \simeq K$. This can be viewed as a special case of the Grothendieck construction.

\begin{rem}
The traditional classification of the group extensions (Schreier theory) requires choosing a set-theoretical section $\Gamma \to G$. This section serves as a `coordinate' to extract the defining data of the monoidal functor $\Gamma \to \Aut(\rmB K)$. In what follows, we adopt a coordinate-free language to perform computations.
\end{rem}

\begin{rem}
Given a monoidal functor $\Gamma \to \Aut(\rmB K)$, the corresponding group extension is the homotopy pullback
\[
\xymatrix{
K \ar[r] \ar@{=}[d] & G \ar[r]^-{s} \ar[d]^{\mathrm{Ad}} \ar@{}[dr]|>>{\ulcorner} & \Gamma \ar[d] \\
K \ar[r]^-{\mathrm{Ad}} & \Aut(K) \ar[r] & \Aut(\rmB K)
}
\]
where both rows are fiber sequences.
\end{rem}

Given a $\Gamma$-graded group $(G,s)$ such that $s$ is surjective and $\ker(s) = K$, now we explicitly write down the corresponding monoidal functor $\Gamma \to \Aut(\rmB K) \simeq \Aut_\Cb^{E_3}(\rep_\Cb(K))$. First, the group $G$ acts on $\Gamma$ via $s$. So we have the action groupoid (homotopy quotient) $\Gamma \dquotient G$:
\bit
\item The objects are elements in $\Gamma$.
\item For any $\sigma,\tau \in \Gamma$, the morphisms $g \colon \sigma \to \tau$ are elements $g \in G$ such that $s(g) \cdot \sigma = \tau$.
\eit
Since $s$ is surjective, $\Gamma \dquotient G$ is connected and thus equivalent to $\rmB K$. There is a right translation $\Gamma$-action on $\Gamma \dquotient G$:
\begin{align*}
\Gamma \times \Gamma \dquotient G & \to \Gamma \dquotient G \\
(\sigma,\tau) & \mapsto \tau \sigma^{-1} .
\end{align*}
This $\Gamma$-action on $\Gamma \dquotient G \simeq \rmB K$ is the one corresponding to $(G,s)$.

For any $(V,\rho) \in \rep_\Cb(K)$, its induced representation $\Cb[G] \otimes_{\Cb[K]} V$ admits a canonical left $G$-action and a $G/K \simeq \Gamma$-grading:
\[
(\Cb[G] \otimes_{\Cb[K]} V)_\sigma \coloneqq \operatorname{span}_\Cb\{g \otimes_{\Cb[K]} v \mid g \in G , \, v \in V , \, s(g) = \sigma \} , \quad \forall \sigma \in \Gamma .
\]
Moreover, the $G$-action and the $\Gamma$-grading are compatible, in the sense that
\[
g \triangleright (\Cb[G] \otimes_{\Cb[K]} V)_\sigma \subseteq (\Cb[G] \otimes_{\Cb[K]} V)_{s(g) \cdot \sigma} .
\]
In other words, $\Cb[G] \otimes_{\Cb[K]} V$ is an object in the equivariantization category $\vect_\Cb[\Gamma]^G$ consisting of finite-dimensional $G$-equivariant $\Gamma$-graded vector spaces, where $G$ acts on $\vect_\Cb[\Gamma]$ via $s \colon G \to \Gamma$. This category is also equivalent to the representation category $\rep_\Cb(\Gamma \dquotient G)$. The assignment $(V,\rho) \mapsto \Cb[G] \otimes_{\Cb[K]} V$ defines a functor $\rep_\Cb(K) \to \vect_\Cb[\Gamma]^G$. It is an equivalence, and the quasi-inverse is given by taking the trivial component of $G$-equivariant $\Gamma$-graded vector spaces. Moreover, the componentwise tensor product defines a symmetric monoidal structure on $\vect_\Cb[\Gamma]^G$ and $\rep_\Cb(K) \simeq \vect_\Cb[\Gamma]^G$ is a symmetric monoidal equivalence. For any $\sigma \in \Gamma$, let $\forget_\sigma \colon \vect_\Cb[\Gamma]^G \to \vect_\Cb$ be the functor of taking the $\sigma$-th component. This is a fiber functor. Let $\omega_\sigma \in \SF_K$ be the composite functor
\[
\rep_\Cb(K) \simeq \vect_\Cb[\Gamma]^G \xrightarrow{\forget_\sigma} \vect_\Cb .
\]
In other words, $\omega_\sigma(V,\rho) \coloneqq (\Cb[G] \otimes_{\Cb[K]} V)_\sigma$. For any $g \in G$, there is a monoidal natural isomorphism $\omega_g \colon \omega_\sigma \Rightarrow \omega_{s(g) \sigma}$ defined by
\begin{align*}
(\omega_g)_{(V,\rho)} \colon (\Cb[G] \otimes_{\Cb[K]} V)_\sigma & \to (\Cb[G] \otimes_{\Cb[K]} V)_{s(g) \sigma} \\
h \otimes_{\Cb[K]} v & \mapsto (gh) \otimes_{\Cb[K]} v .
\end{align*}
Then $\{\omega_\sigma\}_{\sigma \in \Gamma}$ and $\{\omega_g\}_{g \in G}$ define a functor $\omega \colon \Gamma \dquotient G \to \SF_K$, which is clearly an equivalence.

\begin{rem}
Geometrically, $\vect_\Cb[\Gamma]^G$ can be viewed as the category of finite-dimensional $G$-equivariant vector bundles over $\Gamma$. Then $\forget_\sigma$ is literally a `fiber functor': it takes the fiber space over $\sigma \in \Gamma$. The projection $G \to G/K$ is a $G$-equivariant principal $K$-bundle. For any $(V,\rho) \in \rep_\Cb(K)$, its associated vector bundle is defined by
\begin{align*}
G \times_K V = G \times V / \sim & \to G/K \simeq \Gamma \\
[(g,v)] & \mapsto gK \mapsto s(g) 
\end{align*}
where the equivalence relation $\sim$ is generated by $(gk,v) \sim (g,\rho(k)(v))$ for all $k \in K$. This is a $G$-equivariant vector bundle over $\Gamma$. So $(V,\rho) \mapsto G \times_K V$ defines a functor $\rep_\Cb(K) \to \vect_\Cb[\Gamma]^G$. Note that the total space of $G \times_K V$ is exactly the induced representation $\Cb[G] \otimes_{\Cb[K]} V$. Therefore, $\omega_\sigma(V,\rho)$ is exactly the fiber space of $G \times_K V$ over $\sigma \in \Gamma$.
\end{rem}

So we have an equivalence
\[
\rep_\Cb(K) \xrightarrow{\ev} \fun(\SF_K,\vect_\Cb) \xrightarrow{- \circ \omega} \fun(\Gamma \dquotient G,\vect_\Cb) .
\]
Its quasi-inverse is given by
\[
\fun(\Gamma \dquotient G,\vect_\Cb) \xrightarrow{- \circ \iota} \fun(\rmB K,\vect_\Cb) \simeq \rep_\Cb(K) ,
\]
where $\iota \colon \rmB K \to \Gamma \dquotient G$ is the equivalence sending the unique object to the unit of $\Gamma$. So this functor is the evaluation at the unit of $\Gamma$. Then we can explicitly write down the $\Gamma$-action on $\rep_\Cb(K)$ induced by the right translation action on $\Gamma \dquotient G$. This action is a monoidal functor $T \colon \Gamma \to \Aut_\Cb^{E_3}(\rep_\Cb(K))$ defined as follows.
\bit
\item For any $\sigma \in \Gamma$, the symmetric monoidal equivalence $T_\sigma \in \Aut_\Cb^{E_3}(\rep_\Cb(K))$ is the composite functor
\[
\xymatrix{
\rep_\Cb(K) \ar[r]^-{\ev} \ar@{-->}[d]_{T_\sigma} & \fun(\SF_K,\vect_\Cb) \ar[r]^-{- \circ \iota} & \fun(\Gamma \dquotient G,\vect_\Cb) \ar[d]^{- \circ (\sigma^{-1} \triangleright -)} \\
\rep_\Cb(K) & \fun(\rmB K,\vect_\Cb) \ar@{<->}[l]_-{\simeq} & \fun(\Gamma \dquotient G,\vect_\Cb) \ar[l]_-{- \circ \iota}
}
\]
Explicitly, $T_\sigma(V,\rho) \in \rep_\Cb(K)$ is the vector space $\omega_\sigma(V,\rho)$ equipped with the $K$-action induced by $\omega_k$ for all $k \in K$.
\item For $\sigma,\tau \in \Gamma$, the monoidal natural isomorphism $T_\sigma T_\tau \Rightarrow T_{\sigma \tau}$ is defined by
\begin{align*}
T_\sigma(T_\tau(V,\rho)) & \to T_{\sigma \tau}(V,\rho) \\
g \otimes_{\Cb[K]} (h \otimes_{\Cb[K]} v) & \mapsto (gh) \otimes_{\Cb[K]} v .
\end{align*}
It is easy to check that the above map is $K$-equivariant.
\eit

The following result is well-known.

\begin{lem} \label{lem_extension_action_equivariantization_bosonic_linear}
Given the above $\Gamma$-action on $\rep_\Cb(K)$, the equivariantization $\rep_\Cb(K)^\Gamma$ is equivalent to $\rep_\Cb(G)$.
\end{lem}

\pf
The direct sum $\bigoplus_{\sigma \in \Gamma} T_\sigma$ is a monad on $\rep_\Cb(K)$. It maps a $K$-representation $(V,\rho) \in \rep_\Cb(K)$ to
\[
\bigoplus_{\sigma \in \Gamma} T_\sigma(V,\rho) = \Cb[G] \otimes_{\Cb[K]} V ,
\]
which is the induced $G$-representation. Therefore, this monad is induced by the restriction-induction adjunction. The equivariantization category $\rep_\Cb(K)^\Gamma$ is equivalent to the Eilenberg-Moore category of this monad, that is, $\rep_\Cb(G)$.
\epf

\begin{rem}
Lemma \ref{lem_extension_action_equivariantization_bosonic_linear} is a special case of the following well-known result. Let $\CC$ be a category equipped with a left action of a group $G$. For any normal subgroup $N$ of $G$, the equivariantization $\CC^N$ admits a $G/N$-action and $\CC^G \simeq (\CC^N)^{G/N}$. 
\end{rem}

In addition, suppose $z \in K$ is a central element satisfying $z^2 = e$. This element $z$ defines a $\rmB \Zb_2$-action on $\rmB K \simeq \Gamma \dquotient G$, that is, a natural automorphism $\alpha^z$ of the identity functor on $\Gamma \dquotient G$ defined by
\[
\alpha^z_\sigma \coloneqq \psi(\sigma) \cdot z \cdot \psi(\sigma)^{-1} , \quad \sigma \in \Gamma \dquotient G ,
\]
where $\psi \colon \Gamma \to G$ is an arbitrary set-theoretical section. Clearly the $\Gamma$-action on $\Gamma \dquotient G \simeq \rmB K$ is $\rmB \Zb_2$-equivariant if and only if $\alpha^z_\sigma = z$ for all $\sigma \in \Gamma \dquotient G$, which is further equivalent to that $z \in G$ is central by the naturality of $\alpha^z$. Therefore, the extension $(G,z)$ is a super group if and only if the corresponding $\Gamma$-action on $\rmB K$ is $\rmB \Zb_2$-equivariant. By Proposition \ref{prop_auto_srep_G_z_linear}, this super group extension induces a $\Gamma$-action on $\srep_\Cb(K,z)$.

\begin{lem}
The equivariantization $\srep_\Cb(K,z)^\Gamma$ is equivalent to $\srep_\Cb(G,z)$.
\end{lem}

\pf
By Lemma \ref{lem_extension_action_equivariantization_bosonic_linear}, $\srep_\Cb(K,z)^\Gamma$ and $\srep_\Cb(G,z)$ are equivalent as monoidal categories. The statement follows from the fact that the restriction functor $\srep_\Cb(G,z) \to \srep_\Cb(K,z)$ is symmetric.
\epf

\subsection{Classification of symmetric fusion categories over \texorpdfstring{$\Rb$}{R}}

Let $\CC$ be a symmetric fusion category over $\Rb$. If $\Omega \CC \simeq \Cb$, by Lemma \ref{lem_braided_fusion_category_over_Omega} it is a symmetric fusion category over $\Cb$. By Deligne's theorem (see Theorem \ref{thm_Deligne_classificaion_symmetric_fusion_category}) it is equivalent to $\srep_\Cb(G,z) \simeq \srep_{\Cb/\Rb}((G,z) \times \Zb_2^T)$ for some finite super group $(G,z)$.

Now suppose $\Omega \CC \simeq \Rb$. By Theorem \ref{thm_Galois_descent_category}, it is equivalent to the equivariantization of its complexification $\CC_\Cb$ equipped with some semi-linear $\Zb_2$-action. Clearly $\CC_\Cb$ is a symmetric fusion category over $\Cb$, which by Deligne's theorem is classified by super groups.

First we consider the bosonic case, i.e., $\CC_\Cb \simeq \rep_\Cb(K)$ for some finite group $K$. By definition, a symmetric semi-linear $\Zb_2$-action on $\rep_\Cb(K)$ is a monoidal functor $\tilde T \colon \Zb_2 \to \Aut_{\Cb/\Rb}^{E_3}(\rep_\Cb(K))$ such that the following diagram commutes:
\[
\xymatrix{
\Zb_2 \ar[rr]^-{\tilde T} \ar@{=}[dr] & & \Aut_{\Cb/\Rb}^{E_3}(\rep_\Cb(K)) \ar[dl] \\
 & \Zb_2
}
\]
By Lemma \ref{lem_auto_rep_G_semi_linear}, such a semi-linear $\Zb_2$-action on $\rep_\Cb(K)$ is equivalent to a monoidal functor $T \colon \Zb_2 \to \Aut_\Cb^{E_3}(\rep_\Cb(K))$, and
\[
\tilde T_\sigma = T_\sigma \circ \overline{(\cdot)}^\sigma , \quad \sigma \in \Zb_2 .
\]
By the discussion in Section \ref{sec_extension_equivariantization}, the monoidal functor $T$ is induced by a group extension
\[
\xymatrix{
1 \ar[r] & K \ar[r] & G \ar[r]^-{s} & \Zb_2 \ar[r] & 1
}
\]
or equivalently, a $\Zb_2$-graded group $(G,s)$ such that $s$ is surjective and $\ker(s) \simeq K$.

\begin{prop} \label{prop_extension_action_equivariantization_bosonic_semi_linear}
Given a semi-linear $\Zb_2$-action on $\rep_\Cb(K)$ induced by a $\Zb_2$-graded group $(G,s)$ with nontrivial $s$, the equivariantization $\rep_\Cb(K)^{\Zb_2}$ is equivalent to $\rep_{\Cb/\Rb}(G,s)$.
\end{prop}

\pf
The proof is similar to that of Lemma \ref{lem_extension_action_equivariantization_bosonic_linear}. The monad $\bigoplus_{\sigma \in \Zb_2} \tilde T_\sigma$ maps a $K$-representation $(V,\rho) \in \rep_\Cb(K)$ to
\[
\bigoplus_{\sigma \in \Zb_2} \tilde T_\sigma(V,\rho) = \Cb[G,s] \otimes_{\Cb[K]} V .
\]
This space is the same as the induced representation $\Cb[G] \otimes_{\Cb[K]} V$, but the $G$-action is semi-linear because $\tilde T_1 = T_1 \circ \overline{(\cdot)}$ is anti-linear. So this monad is induced by the restriction-induction adjunction between $\rep_\Cb(K) \simeq \LMod_{\Cb[K]}(\vect_\Cb)$ and $\rep_{\Cb/\Rb}(G,s) \simeq \LMod_{\Cb[G,s]}(\vect_\Cb)$.
\epf

For the fermionic case, i.e., $\CC_\Cb \simeq \srep_\Cb(K,z)$ for some finite super group $(K,z)$, we also have a similar result.

\begin{prop} \label{prop_extension_action_equivariantization_fermionic_semi_linear}
Given a semi-linear $\Zb_2$-action on $\srep_\Cb(K,z)$ induced by a $\Zb_2$-graded super group $(G,z,s)$ with nontrivial $s$, the equivariantization $\srep_\Cb(K,z)^{\Zb_2}$ is equivalent to $\srep_{\Cb/\Rb}(G,z,s)$.
\end{prop}

\begin{cor} \label{cor_auto_srep_R_linear}
Let $(G,z,s)$ be a $\Zb_2$-graded finite super group such that $s$ is surjective. Denote $K \coloneqq \ker(s)$. The 2-group $\Aut^{E_3}_\Rb(\srep_{\Cb/\Rb}(G,z,s))$ is equivalent to $\Aut_{\Zb_2 \times \rmB \Zb_2}(\rmB K)$.
\end{cor}

\pf
Combining Proposition \ref{prop_extension_action_equivariantization_fermionic_semi_linear} and Theorem \ref{thm_Galois_descent_category}, we have
\[
\Aut^{E_3}_\Rb(\srep_{\Cb/\Rb}(G,z,s)) \simeq \Aut^{E_3}_\Rb(\srep_\Cb(K,z)^{\Zb_2}) \simeq \Aut^{E_3}_\Cb(\srep_\Cb(K,z))^{\Zb_2} .
\]
By Proposition \ref{prop_auto_srep_G_z_linear}, it is equivalent to $\Aut_{\rmB \Zb_2}(\rmB K)^{\Zb_2} \simeq \Aut_{\Zb_2 \times \rmB \Zb_2}(\rmB K)$.
\epf

\begin{rem} \label{rem_auto_srep_R_linear}
Let $(G,z)$ be a finite super group. Also by Theorem \ref{thm_Galois_descent_category}, we have
\[
\Aut^{E_3}_\Rb(\srep_\Cb(G,z)) \simeq \Aut^{E_3}_\Rb(\vect_\Cb \boxtimes_\Rb \srep_\Cb(G,z))^{\Zb_2} \simeq \Aut^{E_3}_\Rb(\srep_\Cb(G,z) \oplus \srep_\Cb(G,z))^{\Zb_2} ,
\]
where the $\Zb_2$-action on $\srep_\Cb(G,z) \oplus \srep_\Cb(G,z)$ permutes two components. By Proposition \ref{prop_auto_srep_G_z_linear}, it is equivalent to
\[
\Aut_{\rmB \Zb_2}(\rmB G \sqcup \rmB G)^{\Zb_2} \simeq \Aut_{\Zb_2 \times \rmB \Zb_2}(\Zb_2 \times \rmB G) \simeq \Aut_{\Zb_2}(\Zb_2) \times \Aut_{\rmB \Zb_2}(\rmB G) \simeq \Zb_2 \times \Aut_{\rmB \Zb_2}(\rmB G) .
\]
In particular, $\Aut^{E_3}_\Rb(\vect_\Cb) \simeq \Zb_2$ and $\Aut^{E_3}_\Rb(\svect_\Cb) \simeq \Zb_2 \times \rmB \Zb_2$.
\end{rem}

\begin{rem}
Let $\vect_{\Cb/\Rb}$ be the category whose objects are finite-dimensional $\Cb$-vector spaces and morphisms are $\Cb$-linear or anti-linear maps. There is an obvious functor $\vect_{\Cb/\Rb} \to \rmB \Zb_2$ sending $\Cb$-linear maps to $0$ and anti-linear maps to $1$. Then a semi-linear representation $(V,\rho) \in \rep_{\Cb/\Rb}(G,s)$ is the same as a functor $\rho \colon \rmB G \to \vect_{\Cb/\Rb}$ such that the following diagram commutes:
\[
\xymatrix{
\rmB G \ar[rr]^-{\rho} \ar[dr]_{\rmB s} & & \vect_{\Cb/\Rb} \ar[dl] \\
 & \rmB \Zb_2
}
\]
We say that $\rho \colon \rmB G \to \vect_{\Cb/\Rb}$ is a \emph{functor over $\rmB \Zb_2$}.

One can check that $\vect_{\Cb/\Rb}$ is equivalent to the homotopy quotient (see for example \cite[Section 3.2]{HXZ24}) $\vect_\Cb \dquotient \Zb_2$, where the $\Zb_2$-action on $\vect_\Cb$ is given by complex conjugation. Also, $\rmB G$ is equivalent to $\rmB K \dquotient \Zb_2$. By taking the homotopy fiber, a functor $\rmB K \dquotient \Zb_2 \to \vect_\Cb \dquotient \Zb_2$ is equivalent to a $\Zb_2$-equivariant functor $\rmB K \to \vect_\Cb$. Then we obtain another proof of Proposition \ref{prop_extension_action_equivariantization_bosonic_semi_linear}:
\[
\rep_{\Cb/\Rb}(G,s) \simeq \fun_{\Zb_2}(\rmB K,\vect_\Cb) \simeq \fun(\rmB K,\vect_\Cb)^{\Zb_2} \simeq \rep_\Cb(K)^{\Zb_2} .
\]
Similarly, by considering $\svect_{\Cb/\Rb} \simeq \svect_\Cb \dquotient \Zb_2$ we obtain another proof of Proposition \ref{prop_extension_action_equivariantization_fermionic_semi_linear}.
\end{rem}

Combining Proposition \ref{prop_extension_action_equivariantization_bosonic_semi_linear} and \ref{prop_extension_action_equivariantization_fermionic_semi_linear}, we obtain the following classification result.

\begin{thm} \label{thm_classification_symmetric_fusion_R}
Every symmetric fusion category over $\Rb$ is equivalent to $\srep_{\Cb/\Rb}(G,z,s)$ for some $\Zb_2$-graded finite super group $(G,z,s)$.
\end{thm}

\subsection{Fiber functors}

Let $\CC$ be a symmetric fusion category over $\Rb$.

\begin{lem} \label{lem_fiber_functor_over_R_unique}
If $\Omega \CC \simeq \Rb$, there exists a unique $\Rb$-linear fiber functor $\CC \to \svect_\Cb$ up to monoidal natural isomorphism. If $\Omega \CC \simeq \Cb$, there exist exactly two $\Rb$-linear fiber functors $\CC \to \svect_\Cb$ up to monoidal natural isomorphism.
\end{lem}

\pf
By the universal property of complexification, an $\Rb$-linear fiber functor $\CC \to \svect_\Cb$ can be uniquely lifted to a $\Cb$-linear fiber functor $\CC_\Cb \to \svect_\Cb$ (up to monoidal natural isomorphism).

If $\Omega \CC \simeq \Rb$, the complexification $\CC_\Cb$ is a symmetric fusion category over $\Cb$. If $\Omega \CC \simeq \Cb$, the complexification $\CC_\Cb \simeq \vect_\Cb \boxtimes_\Rb \CC$ is the direct sum of two symmetric fusion categories because $\Omega \CC_\Cb \simeq \Cb \otimes_\Rb \Cb \simeq \Cb \oplus \Cb$ as $\Rb$-algebras. Then the statement follows from Deligne's theorem (see Theorem \ref{thm_Deligne_classificaion_symmetric_fusion_category} and \ref{thm_uniqueness_fiber_functor}).
\epf

Let $\tilde \SF_\CC$ be the groupoid defined as follows:
\bit
\item The objects are $\Rb$-linear fiber functors $\CC \to \svect_\Cb$.
\item The morphisms are monoidal natural isomorphisms.
\eit
It admits an action of $\Aut^{E_3}_\Rb(\svect_\Cb) \simeq \Zb_2 \times \rmB \Zb_2$ (see Remark \ref{rem_auto_srep_R_linear}), where the action of the generator of $\Zb_2$ is given by composing with the complex conjugation functor on $\svect_\Cb$.

Define another groupoid $\SF_\CC$ as follows:
\bit
\item The objects are $\Rb$-linear fiber functors $\CC \to \svect_\Cb$.
\item Given two $\Rb$-linear fiber functors $F,F' \colon \CC \to \svect_\Cb$, the morphisms $F \to F'$ are monoidal natural isomorphisms $\forget \circ F \Rightarrow \forget \circ F'$, where $\forget \colon \svect_\Cb \to \svect_\Rb$ is the forgetful functor.
\eit
It admits an action of $\Aut^{E_3}_\Rb(\svect_\Rb) \simeq \rmB \Zb_2$.

\begin{lem} \label{lem_monoidal_natural_isomorphism_semi-linear}
Let $\alpha \colon F \to F'$ be a morphism in $\SF_\CC$. Then it is either $\Cb$-linear or anti-linear, in the sense that $\alpha_x \colon F(x) \to F'(x)$ is $\Cb$-linear or anti-linear for all $x \in \CC$. In particular, there is a functor $\SF_\CC \to \rmB \Zb_2$.
\end{lem}

\pf
The morphism
\[
\kappa \coloneqq \bigl( \Cb \simeq F(\one) \xrightarrow{\alpha_\one} F'(\one) \simeq \Cb \bigr)
\]
is an $\Rb$-linear map. Since $\alpha$ is monoidal, we have the following commutative diagram:
\[
\xymatrix{
\Cb \otimes_\Rb \Cb \ar[r]^-{\simeq} \ar[d]_{\kappa \otimes_\Rb \kappa} & F(\one) \otimes_\Rb F(\one) \ar[r] \ar[d]_{\alpha_\one \otimes_\Rb \alpha_\one} & F(\one) \ar[r]^-{\simeq} \ar[d]^{\alpha_\one} & \Cb \ar[d]^{\kappa} \\
\Cb \otimes_\Rb \Cb \ar[r]^-{\simeq} & F(\one) \otimes_\Rb F(\one) \ar[r] & F(\one) \ar[r]^-{\simeq} & \Cb
}
\]
The composites of the two horizontal arrows are the multiplication of $\Cb$. It follows that $\kappa \colon \Cb \to \Cb$ is an $\Rb$-linear algebra homomorphism. Therefore, $\kappa \in \Gal(\Cb/\Rb) \simeq \Zb_2$ is either identity or the complex conjugation.

Also, for any $x \in \CC$ we have the following commutative diagram:
\[
\xymatrix{
\Cb \otimes_\Rb F(x) \ar[r]^-{\simeq} \ar[d]_{\kappa \otimes_\Rb \alpha_x} & F(\one) \otimes_\Rb F(x) \ar[r] \ar[d]_{\alpha_\one \otimes_\Rb \alpha_x} & F(x) \ar[d]^{\alpha_x} \\
\Cb \otimes_\Rb F(x) \ar[r]^-{\simeq} & F(\one) \otimes_\Rb F(x) \ar[r] & F(x)
}
\]
The composites of the two horizontal arrows are the scalar multiplication of the $\Cb$-vector space $F(x)$. So $\alpha_x$ is $\Cb$-linear if $\kappa$ is identity, and $\alpha_x$ is anti-linear if $\kappa$ is the complex conjugation.

It is clear that $\alpha \mapsto \kappa$ defines a functor $\SF_\CC \to \rmB \Zb_2$.
\epf

There is an obvious $\Zb_2$-equivariant functor $\tilde \SF_\CC \to \SF_\CC$. It induces a functor $\tilde \SF_\CC \dquotient \Zb_2 \to \SF_\CC$, which is clearly $\rmB \Zb_2$-equivariant.

\begin{prop} \label{prop_fiber_functor_groupoid_quotient}
The canonical functor $\tilde \SF_\CC \dquotient \Zb_2 \to \SF_\CC$ is an equivalence such that the following diagram commutes:
\[
\xymatrix{
\tilde \SF_\CC \dquotient \Zb_2 \ar[rr]^-{\simeq} \ar[dr] & & \SF_\CC \ar[dl] \\
 & \rmB \Zb_2
}
\]
where the left arrow is induced by the trivial functor $\tilde \SF_\CC \to \ast$ and the right arrow is defined in Lemma \ref{lem_monoidal_natural_isomorphism_semi-linear}.
\end{prop}

\pf
By definition, the objects of $\tilde \SF_\CC \dquotient \Zb_2$ are $\Rb$-linear fiber functors $\CC \to \svect_\Cb$, and the morphisms are pairs $(\sigma,\tilde \alpha) \colon F \to F'$, where $\sigma \in \Zb_2$ and $\tilde \alpha \colon \overline{(\cdot)}^{\sigma} \circ F \Rightarrow F'$ is a monoidal natural isomorphism. The functor $\tilde \SF_\CC \dquotient \Zb_2 \to \SF_\CC$ acts as identity on objects and sends a morphism $(\sigma,\tilde \alpha)$ to $\forget \circ \tilde \alpha \colon \forget \circ F \Rightarrow \forget \circ F'$. It suffices to show that, for any monoidal natural isomorphism $\alpha \colon \forget \circ F \Rightarrow \forget \circ F'$, there exists a unique morphism $(\sigma,\tilde \alpha) \colon F \to F'$ in $\tilde \SF_\CC \dquotient \Zb_2$ such that $\forget \circ \tilde \alpha = \alpha$.

First we show the uniqueness. Suppose there is such a morphism $(\sigma,\tilde \alpha)$. The condition $\forget \circ \tilde \alpha = \alpha$ implies that $\tilde \alpha_x = \alpha_x$ as maps for all $x \in \CC$. Then by definition, $\sigma$ is $0$ or $1$ if $\alpha_x$ is $\Cb$-linear or anti-linear for all $x \in \CC$, respectively.

Then we show the existence. By Lemma \ref{lem_monoidal_natural_isomorphism_semi-linear}, $\alpha$ is either $\Cb$-linear or anti-linear. If $\alpha$ is $\Cb$-linear, it already defines a monoidal natural isomorphism $\alpha \colon F \Rightarrow F'$. So $(0,\alpha) \colon F \to F'$ is a morphism in $\tilde \SF_\CC \dquotient \Zb_2$. If $\alpha$ is anti-linear, it can be viewed as a monoidal natural isomorphism $\alpha \colon \overline{(\cdot)} \circ F \Rightarrow F'$. So $(1,\alpha) \colon F \to F'$ is a morphism in $\tilde \SF_\CC \dquotient \Zb_2$.

The commutativity of the above diagram is obvious.
\epf

Let us compute the super groupoids $\tilde \SF_\CC$ and $\SF_\CC$ for $\CC = \srep_{\Cb/\Rb}(G,z,s)$, where $(G,z,s)$ is a $\Zb_2$-graded finite super group.

Suppose $s$ is trivial, i.e., $\CC \simeq \srep_\Cb(G,z)$. Then by Lemma \ref{lem_monoidal_natural_isomorphism_semi-linear}, $\tilde \SF_\CC$ is equivalent to $\Zb_2 \times \rmB G$, on which the $\rmB \Zb_2$-action is given by $z$. So $\SF_\CC \simeq \tilde \SF_\CC \dquotient \Zb_2 \simeq \rmB G$.

Now suppose $s$ is nontrivial. By Lemma \ref{lem_fiber_functor_over_R_unique}, $\tilde \SF_\CC$ is connected. Then by Proposition \ref{prop_fiber_functor_groupoid_quotient}, $\SF_\CC$ is also connected. So we only need to compute the automorphism groups of the forgetful functors $\omega \colon \srep_{\Cb/\Rb}(G,z,s) \to \svect_\Cb$ and $\forget \circ \omega \colon \srep_{\Cb/\Rb}(G,z,s) \to \svect_\Cb \to \svect_\Rb$.

For any finite-dimensional $\bk$-algebras $A,B$ and $\bk$-algebra homomorphism $\phi \colon A \to B$, it is well-known that the endomorphism algebra of the restriction functor $\LMod_B(\vect_\bk) \to \LMod_A(\vect_\bk)$ is the centralizer
\[
Z(\phi) \coloneqq \{z \in B \mid \phi(a) z = z \phi(a) , \, \forall a \in A \} .
\]
So every natural endomorphism of $\omega$ must be $\alpha^x$ for some $x \in Z(\Cb[\Zb_2] \xrightarrow{z} \Cb[G,s]) = \Cb[\ker(s)]$, where
\[
\alpha^x_{(V,\rho)} \coloneqq \rho(x) \colon V \to V .
\]
It is monoidal if $\alpha^x_{(V,\rho)} \otimes_\Cb \alpha^x_{(W,\sigma)} = \alpha^x_{(V,\rho) \otimes_\Cb (W,\sigma)}$ for all $(V,\rho),(W,\sigma) \in \srep_{\Cb/\Rb}(G,z,s)$. Taking $(V,\rho)$ and $(W,\sigma)$ to be the regular semi-linear representation $\Cb[G,s]$ (or any faithful semi-linear representation), we obtain
\[
x \otimes_\Cb x = \Delta(x) ,
\]
where the map $\Delta \colon \Cb[G,s] \to \Cb[G,s] \otimes_\Cb \Cb[G,s]$ sends $g \in G$ to $g \otimes_\Cb g$. Therefore, $\alpha^x$ is a monoidal natural transformation if and only if $x \in \ker(s)$. So $\tilde \SF_\CC \simeq \rmB \ker(s)$.

Similarly, every natural endomorphism of $\forget \circ \omega$ must be $\alpha^x$ for some $x \in Z(\Rb[\Zb_2] \xrightarrow{z} \Cb[G,s]) = \Cb[G,s]$, which is monoidal if and only if $x \in G$. So $\SF_\CC \simeq \rmB G$. Clearly the composite functor $\rmB G \hookrightarrow \SF_\CC \to \rmB \Zb_2$ is given by $\rmB s$. Then by Proposition \ref{prop_fiber_functor_groupoid_quotient}, the $\Zb_2$-action on $\tilde \SF_\CC \simeq \rmB \ker(s)$ is the one corresponding to the group extension $1 \to \ker(s) \to G \xrightarrow{s} \Zb_2 \to 1$. Moreover, the $\rmB \Zb_2$-action on $\tilde \SF_\CC \simeq \rmB \ker(s)$ and $\SF_\CC \simeq \rmB G$ are both given by $z \in G$.

These computations together with Theorem \ref{thm_classification_symmetric_fusion_R} gives the following result.

\begin{thm} \label{thm_symmetric_fusion_R_fiber_functor_groupoid}
Let $\CC$ be a symmetric fusion category over $\Rb$. The canonical functor $\CC \to \fun_{\Zb_2 \times \rmB \Zb_2}(\tilde \SF_\CC,\svect_\Cb)$ is an equivalence.
\end{thm}

We denote the 2-groupoid of $\Rb$-linear symmetric multi-fusion categories, $\Rb$-linear symmetric monoidal equivalences and monoidal natural isomorphisms by $\Alg_{E_3}^\rig(2\vect_\Rb)^\times$, and denote the 2-groupoid of finite groupoids equipped with a $\Zb_2 \times \rmB \Zb_2$-action, equivariant equivalences and equivariant natural isomorphisms by $\LMod_{\Zb_2 \times \rmB \Zb_2}(\mathrm{fGrpd})^\times$.

\begin{thm} \label{thm_symmetric_fusion_category_R_super_groupoid_Galois_action}
The functors
\begin{align*}
(\Alg_{E_3}^\rig(2\vect_\Rb)^\times)^\op & \to \LMod_{\Zb_2 \times \rmB \Zb_2}(\mathrm{fGrpd})^\times & (\LMod_{\Zb_2 \times \rmB \Zb_2}(\mathrm{fGrpd})^\times)^\op & \to \Alg_{E_3}^\rig(2\vect_\Rb)^\times \\
\CC & \mapsto \tilde \SF_\CC & X & \mapsto \fun_{\Zb_2 \times \rmB \Zb_2}(X,\svect_\Cb)
\end{align*}
define an opposite equivalence between the 2-groupoids $\Alg_{E_3}^\rig(2\vect_\Rb)^\times$ and $\LMod_{\Zb_2 \times \rmB \Zb_2}(\mathrm{fGrpd})^\times$.
\end{thm}

\pf
By Theorem \ref{thm_classification_symmetric_fusion_R} and Theorem \ref{thm_symmetric_fusion_R_fiber_functor_groupoid}, the 2-functor $X \mapsto \fun_{\Zb_2 \times \rmB \Zb_2}(X,\svect_\Cb)$ is essentially surjective. By Corollary \ref{cor_auto_srep_R_linear} and Remark \ref{rem_auto_srep_R_linear}, this 2-functor is fully faithful. Also by Theorem \ref{thm_symmetric_fusion_R_fiber_functor_groupoid}, the quasi-inverse is given by $\CC \mapsto \tilde \SF_\CC$.
\epf

Finally, we determine which symmetric fusion categories admit fiber functors to $\vect_\Rb$ or $\svect_\Rb$.

\begin{thm}
A symmetric fusion category $\CC$ over $\Rb$ admits an $\Rb$-linear fiber functor to $\vect_\Rb$ (that is, Tannakian over $\Rb$) if and only if it is equivalent to $\rep_{\Cb/\Rb}(K \rtimes \Zb_2^T)$ for some finite group $K$ equipped with a $\Zb_2$-action.
\end{thm}

\pf
If there exists such a fiber functor, it induces an $\Rb$-linear algebra homomorphism $\Omega \CC \to \Omega \vect_\Rb \simeq \Rb$. So $\CC$ admits an $\Rb$-linear fiber functor to $\vect_\Rb$ only if $\Omega \CC \simeq \Rb$.

By Theorem \ref{thm_Galois_descent_category}, an $\Rb$-linear fiber functor $\CC \to \vect_\Rb$ is equivalent to a $\Zb_2$-equivariant $\Cb$-linear fiber functor $\CC_\Cb \to \vect_\Cb$. Assume that the complexification $\CC_\Cb$ is equivalent to $\rep_\Cb(K)$ for some finite group $K$ and the semi-linear $\Zb_2$-action $\tilde T \colon \Zb_2 \to \Aut^{E_3}_{\Cb/\Rb}(\rep_\Cb(K))$ corresponds to a group extension
\[
\xymatrix{
1 \ar[r] & K \ar[r] & G \ar[r]^-{s} & \Zb_2 \ar[r] & 1
}
\]
It suffices to determine whether the forgetful functor $\forget \colon \rep_\Cb(K) \to \vect_\Cb$ is $\Zb_2$-equivariant.

A $\Zb_2$-equivariant structure on $\forget$ consists of natural isomorphisms
\[
\overline{(\cdot)}^{\sigma} \circ \forget \circ \tilde T_\sigma \Rightarrow \forget , \quad \sigma \in \Zb_2 ,
\]
satisfying certain coherence conditions. Note that $F \mapsto \overline{(\cdot)}^{\sigma} \circ \forget \circ \tilde T_\sigma$ defines a semi-linear $\Zb_2$-action on the functor category $\fun_\Cb(\rep_\Cb(K),\vect_\Cb)$, and the Yoneda embedding $\rep_\Cb(K)^\op \to \fun_\Cb(\rep_\Cb(K),\vect_\Cb)$ is $\Zb_2$-equivariant. Therefore, a $\Zb_2$-equivariant structure on $\forget \simeq \rep_\Cb(K)(\Cb[K],-)$ is equivalent to a $\Zb_2$-equivariant structure on the regular representation $\Cb[K] \in \rep_\Cb(K)$. By Proposition \ref{prop_extension_action_equivariantization_bosonic_semi_linear}, it is equivalent to a semi-linear $(G,s)$-action on $\Cb[K]$ such that the the action of elements in $K$ is given by left translation. Moreover, since the monoidal structure of $\forget$ is induced by the diagonal coalgebra structure of $\Cb[K]$, a $\Zb_2$-equivariant structure on $\forget$ is monoidal if and only if the corresponding semi-linear $(G,s)$-action on $\Cb[K]$ preserves the coalgebra structure. Clearly it is the same as a left $G$-action on the set $K$ such that the action of elements in $K$ is given by left translation.

Now suppose $\phi \colon G \times K \to K$ is a left $G$-action on the set $K$ such that $\phi(k_1,k_2) = k_1 k_2$ for all $k_1,k_2 \in K$. Define
\[
H \coloneqq \{g \in G \mid \phi(g,e) = e\} .
\]
Then $H$ is a subgroup of $G$. One can show that:
\bit
\item $K \cap H = \{e\}$: for any $k \in K$, if $k \in H$, then $e = \phi(k,e) = k$;
\item $G = KH$: for any $g \in G$, we have $k \coloneqq \phi(g,e) \in K$ and $k^{-1}g \in H$ because $\phi(k^{-1}g,e) = \phi(k^{-1},\phi(g,e)) = \phi(k^{-1},k) = k^{-1} k = e$.
\eit
Therefore, $G$ is the semidirect product of $K$ and $H \simeq G/K \simeq \Zb_2$. In other words, $(G,s) \simeq K \rtimes \Zb_2^T$.
\epf

Similarly one can prove the following result.

\begin{thm}
A symmetric fusion category $\CC$ over $\Rb$ admits an $\Rb$-linear fiber functor to $\svect_\Rb$ if and only if it is equivalent to $\srep_{\Cb/\Rb}((K,z) \rtimes \Zb_2^T)$ for some finite super group $(K,z)$ equipped with a $\Zb_2$-action.
\end{thm}

\begin{rem}
In general, $\Rb$-linear fiber functors from a symmetric fusion category to $\vect_\Rb$ or $\svect_\Rb$ (if exist) are not unique up to monoidal isomorphism. For example, $\Rb$-linear fiber functor $\rep_\Rb(G) \to \vect_\Rb$ are classified by $G$-torsors over $\Rb$ \cite[Theorem 3.2]{DM82}. When $G = \Zb_2$, since $H^1(\Rb;\Zb_2) \simeq \Zb_2$, there are two $\Rb$-linear fiber functors $\rep_\Rb(\Zb_2) \to \rep_\Rb$ up to monoidal isomorphism. The nontrivial one sends $(V,\rho) \in \rep_\Rb(\Zb_2)$ to $(\Cb \otimes_\Rb V)^{\Zb_2}$, where $\Zb_2$ acts on $\Cb$ by complex conjugation.
\end{rem}

\subsection{Classification over characteristic zero fields}

Most of the results in this paper can be generalized to any field $\bk$ of characteristic zero. Therefore, if $[\overline{\bk} : \bk] < \infty$, we can similarly classify symmetric fusion categories over $\bk$. However, such a field is either algebraically closed or a real closed field. So the classification is essentially the same as $\bk = \Cb$ or $\Rb$.

When $\overline{\bk}/\bk$ is an infinite Galois extension, the classification is more interesting. In this case the Galois group $\Gamma = \Gal(\overline{\bk}/\bk)$ is a profinite group, and we need to require that the $\Gamma$-actions on vector spaces or categories are smooth (continuous). Then symmetric fusion categories over $\bk$ should be classified by smooth semi-linear $\Gamma$-actions on $\srep_{\bar \bk}(G,z)$. We conjecture the classification as follows.

\begin{conj}
Every symmetric fusion category $\CC$ over $\bk$ is equivalent to $\srep_{\Eb/\Fb}(G,z,s)$, where $\Fb/\bk$ is a finite extension, $\Eb/\Fb$ is a finite Galois extension and $(G,z,s)$ is a $\Gal(\Eb/\Fb)$-graded finite super group such that $s$ is surjective. Moreover, the data $(\Eb,\Fb,G,z,s)$ is uniquely determined by $\CC$ up to isomorphism.
\end{conj}

We will study the classification in future works.

\bibliography{Top}

\end{document}